\documentclass{amsart}
\usepackage{amssymb}
\usepackage[all]{xy}
\newcommand\datver[1]{\def\datverp%
 {\par\boxed{\boxed{\text{#1; Run: \today}}}}}

\newcommand\cf{\textit{cf.}}
\newcommand\cho[1]{\beta_{#1}}  
\newcommand\susn[1]{\operatorname{sus}(#1)}
\newcommand\psusn[1]{\operatorname{psus}(#1)}
\newcommand\spe[1]{\circ_{#1}}
\newcommand\fsp{*}
\newcommand\wq[2]{\operatorname{iso}(#1,#2) }
\renewcommand\det{\operatorname{det}}
\newcommand\deta{\operatorname{det_a}}
\newcommand\Ld{\operatorname{L}^{2}} 
\newcommand\br{\rho_{r}}
\newcommand\bs{\rho_{\sus}}
\newcommand\bt{\rho_{\tau}}

\newcommand\pt{\phi-p} 
\newcommand\De{{}^{\epsilon}\hat{D}}
\newcommand\LIM{\operatornamewithlimits{LIM}}

\newcommand\boxb[1]{\square_b}

\newcommand\ff{\operatorname{ff}}

\numberwithin{equation}{section}
\renewcommand{\theequation}{\arabic{equation}}
\let\Appendix\appendix
\newcommand\paperbody%
        {\renewcommand{\theequation}{\thesection.\arabic{equation}}}
\newcommand\firstpart%
        {\renewcommand{\thesection}{A\arabic{section}}
        \renewcommand{\theequation}{\thesection.\arabic{equation}}}
\newcommand\secondpart%
        {\renewcommand{\thesection}{B\arabic{section}}
        \renewcommand{\theequation}{\thesection.\arabic{equation}}}
\renewcommand\appendix{\Appendix%
\renewcommand{\theequation}{A.\arabic{equation}}
\numberwithin{equation}{section}}

\newtheorem{lemma}{Lemma}[section]
\newtheorem{proposition}[lemma]{Proposition}
\newtheorem{corollary}[lemma]{Corollary}
\newtheorem{theorem}{Theorem}
\newtheorem{theorem_int}{Theorem}
\newtheorem{non-theorem}[lemma]{Non-Theorem}

\theoremstyle{remark}

\newtheorem{definition}[lemma]{Definition}

\newcommand\iso{\operatorname{iso}}

\newcommand\coH{\operatorname{H}}

\newcommand\cFTs{{}^{\Phi}\overline{T}\kern-1pt{}^*}

\newcommand\sus{\operatorname{sus}}
\newcommand\psus{\operatorname{psus}}

\newcommand\even{\text{even}}
\newcommand\odd{\text{odd}}

\newcommand\Tr{\operatorname{Tr}}

\newcommand\bTr{\overline{\operatorname{Tr}}}

\newcommand\Det{\operatorname{Det}}
\newcommand\Deta{\operatorname{Det_a}}
\newcommand\tdet{\operatorname{\widetilde{det}}}
\newcommand\adet{\operatorname{{det_a}}}
\newcommand\adn[1]{\operatorname{a(#1)}}

\newcommand\Ma[2]{\operatorname{M}(#1,#2)}
\newcommand\GL[2]{\operatorname{GL}(#1,#2)}

\newcommand\tq{\tilde{q}}

\newcommand\com[1]{\overline{#1}}
\newcommand\pcom[2]{{}^{#1}{\com{#2}}}

\newcommand\ie{i\@.e\@. }

\newcommand\cA{\mathcal{A}}

\newcommand\cH{\mathcal{H}}
\newcommand\cE{\mathcal{E}}

\newcommand\cG{\mathcal{G}}
\newcommand\cI{\mathcal{I}}

\newcommand\cP{\mathcal{P}}

\newcommand\cU{\mathcal{U}}

\newcommand\bbC{\mathbb C}

\newcommand\bbE{\mathbb E}

\newcommand\bbN{\mathbb N}

\newcommand\bbR{\mathbb R}
\newcommand\bbS{\mathbb S}

\newcommand\bbZ{\mathbb Z}

\newcommand\cF{\mathcal F}
\newcommand\cS{\mathcal S}
\newcommand\CIc{{\mathcal{C}}^{\infty}_c}

\newcommand\CI{{\mathcal{C}}^{\infty}}

\newcommand\Diag{\operatorname{Diag}}

\newcommand\cFNs{{}^{\Phi}\overline N\kern-1pt{}^*}

\newcommand\ind{\operatorname{ind}}
\newcommand\res{\operatorname{res}}
\newcommand\wn{\operatorname{wn}}
\newcommand\Fr{\operatorname{Fr}}

\newcommand\reg{\operatorname{reg}}
\newcommand\tr{\operatorname{tr}}

\newcommand\Hom{\operatorname{Hom}}
\newcommand\Aut{\operatorname{Aut}}
\newcommand\EG{\operatorname{E}G}
\newcommand\BG{\operatorname{B}G}
\newcommand\Id{\operatorname{Id}}

\newcommand\ci{${\mathcal{C}}^\infty$}

\newcommand\dCI{\dot{\mathcal{C}}^{\infty}}

\newcommand\pa{\partial}

\newcommand\sgn{\operatorname{sgn}}

\newcommand\cl{\operatorname{cl}}

\renewcommand\Re{\operatorname{Re}}

\newcommand\Mand{\text{ and }}

\newcommand\Mfor{\text{ for }}

\newcommand\Min{\text{ in }}

\newcommand\Mwhere{\text{ where }}

\datver{1.0B; Revised: 15-6-2006}
\begin{document}
\title[Periodicity]
{Perioditicity and the Determinant Bundle}

\author{Richard Melrose}
\address{Department of Mathematics, Massachusetts Institute of Technology}
\email{rbm@math.mit.edu}
\author{Fr\'ed\'eric Rochon}
\address{Department of Mathematics, State University of New York, Stony Brook}
\email{rochon@math.sunysb.edu}
\dedicatory{\datverp}
\thanks{The first author acknowledges the support of the National Science
  Foundation under grant DMS0408993, the second author acknowledges support
  of the Fonds qu\'{e}b\'{e}cois sur la nature et les technologies and
  NSERC while part of this work was conducted}.
\begin{abstract} The infinite matrix `Schwartz' group $G^{-\infty}$ is a
  classifying group for odd K-theory and carries Chern classes in each odd
  dimension, generating the cohomology. These classes are closely related
  to the Fredholm determinant on $G^{-\infty}.$ We show that while the
  higher (even, Schwartz) loop groups of $G^{-\infty},$ again classifying
  for odd K-theory, do \emph{not} carry multiplicative determinants
  generating the first Chern class, `dressed' extensions,
  corresponding to a star product, do carry such functions. We use these to
  discuss Bott periodicity for the determinant bundle and the eta
  invariant. In so doing we relate two distinct extensions of the eta
  invariant, to self-adjoint elliptic operators and to elliptic invertible
  suspended families and show that the corresponding $\tau$ invariant is a
  determinant in this sense.
\end{abstract}
\maketitle

\tableofcontents

\section*{Introduction}

The Fredholm determinant is a character for the group of invertible
operators of the form $\Id+T$ with $T$ of trace class on a Hilbert
space. Transferred to invertible operators of the form $\Id+A$ with $A$
smoothing on the compact fibres of a fibration it induces the determinant
bundle of families of elliptic pseudodifferential operators. For suspended
families of smoothing operators, depending in a Schwartz fashion on an even
number of Euclidean parameters, we introduce an adiabatic
determinant with similar topological properties and use it to prove
periodicity properties for the determinant bundle. The corresponding
suspended eta invariants are also discussed and in a subsequent paper will
be used to describe cobordism of the determinant bundle in a
pseudodifferential setting, extending the result of Dai and Freed
\cite{Dai-Freed1} that the eta invariant in the interior defines a
trivialization of the determinant bundle on the boundary.

The basic notion of determinant is that on finite rank matrices. If $\Ma
N\bbC$ is the algebra of $N\times N$ complex matrices then the determinant
is the entire (polynomial) multiplicative map 
\begin{equation*}
\det:\Ma N\bbC\longrightarrow \bbC,\ \det(AB)=\det(A)\det(B)
\label{perdet.13}\end{equation*}
which is determined by the condition on its derivative at the identity
\begin{equation*}
\frac{d}{ds}\det(\Id+sA)\big|_{s=0}=\Tr(A),\ A\in\Ma N\bbC.
\label{perdet.14}\end{equation*}
It has the fundamental property that $\det(A)\not=0$ is equivalent to the
invertibility of $A,$ so 
\begin{equation*}
\GL N\bbC=\{A\in\Ma N\bbC;\det(A)\not=0\}={\det}^{-1}(\bbC^*).
\label{perdet.15}\end{equation*}

As is well-known, such a map into $\bbC^*$ determines, through the winding
number, an integral 1-cohomology class: 
\begin{equation}
\alpha (c)=\wn(\det:c\longrightarrow \bbC^*),\ \alpha \in\coH^1(\GL N\bbC;\bbZ).
\label{perdet.16}\end{equation}
Conversely for any path-connected space
\begin{equation*}
H^1(X;\bbZ)\equiv\{\alpha:\pi_1(X)\longrightarrow \bbZ;\ \alpha (c_1\circ c_2)
=\alpha (c_1)+\alpha (c_2)\}
\label{fipomb2.1}\end{equation*}
so each integral 1-cohomology class may be represented by a
continuous function $f:X\longrightarrow \bbC^*$ such that $\alpha(c)$ is
the winding number of $f$ restricted to a curve representing $c.$ Even if
$X$ is a group and the class is invariant, it may not be possible to choose
this function to be multiplicative.

Each integral 1-cohomology class on $X$ may also be represented as the
obstruction to the triviality of a principal $\bbZ$ bundle over $X.$ Such a
bundle, with total space $P,$ always admits a `connection' in the sense of
a map $h:P\longrightarrow \bbC$ such that $h(np)=h(p)+n$ for the action of
$n\in\bbZ.$ Given appropriate smoothness, the function on $X$ associated to the
connection, $f=\exp(2\pi ih),$ fixes the obstruction 1-class as a deRham form
\begin{equation*}
\alpha =\frac1{2\pi i}f^{-1}df=dh.
\label{fipomb2.2}\end{equation*}
In particular the triviality of the $\bbZ$-bundle is equivalent to the
existence of a continuous (normalized) logarithm for $f,$ that is a
function $l:X\longrightarrow \bbC$ such that $h-\phi^*l$ is locally
constant, where $\phi:P\longrightarrow X$ is the bundle projection.

Returning to the basic case of the matrix algebra and $\GL N\bbC,$
these spaces can be naturally included in the `infinite matrix algebra'
which we denote abstractly $\Psi^{-\infty}.$ For the moment we identify 
\begin{equation*}
\Psi^{-\infty}=\{a:\bbN^2\longrightarrow
\bbC;\sup_{i,j\in\bbN}(i+j)^k|a_{ij}|<\infty\ \forall\ k\in\bbN\}.
\label{perdet.17}\end{equation*}
The algebra structure is just the extension of standard matrix multiplication 
\begin{equation*}
(ab)_{ij}=\sum\limits_{l=1}^\infty a_{il}b_{lj}.
\label{perdet.18}\end{equation*}
Now, although $\Ma N\bbC\longrightarrow \Psi^{-\infty}$ is included as the
subalgbra with $a_{ij}=0$ for $i,j>N,$ for the determinant this is not
natural, in part because $\Psi^{-\infty}$ is non-unital. Namely, we
consider instead the isomorphic space $\Id+\Psi^{-\infty}$ which may be
identified with $\Psi^{-\infty}$ with the product 
\begin{equation*}
a\circ b=a+b+ab. 
\label{perdet.19}\end{equation*}
Then the inclusion 
\begin{equation*}
\Ma N\bbC\ni a\longmapsto (\Id-\pi_N)+\pi_Na\pi_N
\label{perdet.20}\end{equation*}
is multiplicative and the determinant is consistent for all $N$ with the
Fredholm determinant which is the entire multiplicative function
\begin{equation*}
{\det}_{\Fr}:\Id+\Psi^{-\infty}\longrightarrow \bbC
\label{perdet.21}\end{equation*}
satisfying the normalization
\begin{equation*}
\frac{d}{ds}{\det}_{\Fr}(\Id+sa)\big|_{s=0}=
\Tr(a)=\sum\limits_{i=1}^\infty a_{ii}.
\label{perdet.22}\end{equation*}
Again for $a\in\Psi^{-\infty}$ the condition $\det_{\Fr}(\Id+a)\not=0$ is
equivalent to the existence of an inverse $\Id+b,$ $b\in\Psi^{-\infty}$ and
this defines the topological group 
\begin{equation*}
G^{-\infty}=\{\Id+a;a\in\Psi^{-\infty},\ {\det}_{\Fr}(\Id+a)\not=0\}
\label{perdet.23}\end{equation*}
in which the $\GL N\bbC$ are included as subgroups. Since these determinants
are consistent we generally drop the distinction between the finite and
Fredholm determinants.

Now, $G^{-\infty}$ is a classfying group for odd K-theory, 
\begin{equation*}
K^1(X)=\Pi_0\{f:X\longrightarrow G^{-\infty}\}
\label{perdet.24}\end{equation*}
where the maps can be taken to be either continuous or smooth. As such, 
\begin{equation*}
\Pi_l(G^{-\infty})=\begin{cases}\{0\}& l\text{ even}\\
\bbZ& l\text{ odd}.\end{cases}
\label{perdet.25}\end{equation*}
The odd Chern forms (see for example \cite{fipomb}), 
\begin{equation}
\cho{2k-1}= \frac{1}{(2\pi i)^{k}} \frac{(k-1)!}{(2k-1)!}
\Tr[ ((\Id+a)^{-1}da)^{2k-1} ],\ k\in\bbN,
\label{perdet.89}\end{equation}
give an explicit isomorphism
\begin{equation}
h_{2k-1}:\Pi_{2k-1}(G^{-\infty})\ni[f]\longmapsto \int_{\bbS^{2k-1}} 
f^{*}\cho{2k-1}\in\bbZ
\label{co.1}\end{equation}
where $[f]\in \Pi_{2k-1}(G^{-\infty})$ is represented by a smooth map
$f: \bbS^{2k-1}\longrightarrow G^{-\infty}.$ The cohomology classes
$[\cho{2k-1}]\in H^{2k-1}(G^{-\infty});\bbC)$ generate 
$H^{*}(G^{-\infty};\bbC)$ as 
an \emph{exterior algbra} over $\bbC$ 
\begin{equation*}
H^{*}(G^{-\infty};\bbC) =
\Lambda_{\bbC}( \cho{1}, \cho{3}, \ldots, \cho{2k-1}, \ldots ).
\label{perdet.90}\end{equation*}
However, as noted by Bott and Seeley in \cite{Bott-Seeley}, even though
they give integers 
when integrated over the corresponding spherical homology class, the
 classes $[\cho{2k-1}]$ are not all integral. When
$k=1$, the isomorphism \eqref{co.1} shows that  
as in the case of the matrix groups, a loop along which the winding number
of the determinant is $1$ generates $\Pi_1(G^{-\infty})$ and
$\coH^1(G^{-\infty})$ is generated by the deRham class 
\begin{equation}
\alpha =\frac1{2\pi i}\Tr((\Id+a)^{-1}da).
\label{perdet.27}\end{equation}

Bott periodicity corresponds to the fact that the (reduced) loop groups of
$G^{-\infty}$ are also classifying spaces for odd or even
K-theory. Consistent with the `smooth' structure championed here, we
consider loop groups of `Schwartz' type. In fact we can first identify
$\Psi^{-\infty}$ above as the expansion of an operator with respect to the
eigenvectors of the harmonic oscillator on $\bbR^n$ to identify 
\begin{equation*}
\Psi^{-\infty}\longleftrightarrow \Psi^{-\infty}(\bbR^n)=\mathcal{S}(\bbR^{2n})
\label{perdet.29}\end{equation*}
where the product on $\mathcal{S}(\bbR^{2n})$ is the operator product 
\begin{equation*}
(ab)(x,y)=\int_{\bbR^n}a(x,z)b(z,y)dz.
\label{perdet.30}\end{equation*}
With this identification the loop groups become 
\begin{equation*}
G^{-\infty}_{\susn{p}}(\bbR^n)=\{f:\bbR^p\longrightarrow G^{-\infty}(\bbR^n);
f=\Id+a,\ a\in\mathcal{S}(\bbR^{p+2n})\}.
\label{perdet.28}\end{equation*}
Thus $G^{-\infty}_{\susn{p}}$ is a classifying group for K-theory of the
parity opposite to that of $p.$ In fact we may regard $G^{-\infty}_{\susn{p}}$ as
classifying for the groups $K^{-p-1}$ and Bott periodicity as giving the
identification between these for all even and all odd orders.

The analogues of the forms \eqref{perdet.89} are given by
\begin{equation*}
   \cho{2k-1-p}^{(p)}(f)= \int_{\bbR^{p}} f^{*} \cho{2k-1}, \quad
p\le 2k-1, \; k\in \bbN.
\label{forms.1}\end{equation*}
For $p=1$, this gives the even forms
\begin{equation*}
   \cho{2k}^{(1)}= \frac{1}{(2\pi i)^{k+1}}\frac{k!}{(2k)!}
  \int_{\bbR} \Tr \left[ (a^{-1}da)^{2k} a^{-1}\frac{da}{d\tau} \right]
        d\tau, \quad k\in \bbN_{0},
\label{forms.2}\end{equation*}
where $\tau$ is the suspension parameter (\cf \cite{fipomb}).

It is equally possible to use the eigenbasis of a Laplacian on the sections
of a vector bundle over a compact Riemannian manifold without boundary (or
of any self-adjoint elliptic pseudodifferential operator of positive order)
to identify $\Psi^{-\infty}$ with $\Psi^{-\infty}(X;E),$ the space of
smoothing operators. Then the loop groups are realized as
\begin{multline*}
G^{-\infty}_{\susn{p}}(X;E)=\\
\{f:\bbR^p\longrightarrow G^{-\infty}(X;E);
f=\Id+a,\ a\in\mathcal{S}(\bbR^p\times X\times X;\Hom(E))\}.
\label{perdet.31}\end{multline*}
Here, the space of Schwartz sections is defined for any vector bundle which
is the pull-back to $\bbR^p\times Z$ of a vector bundle over a compact
manifold $Z.$ 

Now, the basic issue considered here is the existence of a determinant on
the spaces $G^{-\infty}_{\susn{2k}}.$ One can simply look for a smooth
multiplicative function which generates the 1-dimensional homology through
the winding number formula \eqref{perdet.16}. In what is really the
opposite side of the `Miracle of the loop group' of Pressley and Segal
\cite{Pressley-Segal1} there is in fact no such function as soon as
$k>0.$ As we show below, there is a multiplicative function closely related
to the determinant but which has a global logarithm (if $k>0).$ However, as
we also show below, there is a determinant function, the `adiabatic determinant'
in this sense provided the group $G^{-\infty}_{\susn{2k}}$ is `dressed' by
replacing it by an extension with respect to a star product, of which
$G^{-\infty}_{\susn{2k}}$ is the principal term. This extension is
homotopically trivial, \ie still gives a classifying space for K-theory.

More precisely, consider the space $\Psi^{k}_{\susn{2n}}(X;E)[[\varepsilon]]$
of formal power series
\begin{equation*}
    \sum_{\mu=0}^{\infty} a_{\mu}\varepsilon^{\mu}, \quad 
                         a_{\mu}\in \Psi^{k}_{\susn{2n}}(X;E)       
\end{equation*}
equipped with the star-product
\begin{equation}
\begin{aligned}
  (A\fsp B)(u)&= \left( \sum_{\mu=0}^{\infty} a_{\mu}\varepsilon^{\mu}\right)
    \fsp \left( \sum_{\nu=0}^{\infty} b_{\nu}\varepsilon^{\nu}\right) \\
   &=  \sum_{\mu=0}^{\infty} \sum_{\nu=0}^{\infty} 
        \varepsilon^{\mu+\nu} 
       \left.\left( \sum_{p=0}^{\infty} \frac{i\varepsilon^{p}}{2^{p}p!} 
           \omega(D_{v},D_{w})^{p} a_{\mu}(v) b_{\nu}(w) \right)
      \right|_{v=w=u}
\end{aligned}
\label{int.1}\end{equation}
for $A, B\in \Psi^{*}_{\susn{2n}}(X;E)[[\varepsilon]]$, where $\omega$ is the 
standard symplectic form on $\bbR^{2n}.$ This gives a corresponding group
\begin{multline*}
G^{-\infty}_{\susn{2n}}(X;E)[[\varepsilon]]= 
\{\Id +Q;\ Q\in \Psi^{-\infty}_{\susn{2n}}(X;E)[[\varepsilon]],\\
\exists\ P\in \Psi^{-\infty}_{\susn{2n}}(X;E)[[\varepsilon]],\
(\Id+Q)\fsp(\Id+ P)= \Id \in \Psi^{0}_{\susn{2n}}(X;E)[[\varepsilon]]\}
\label{int.2}\end{multline*}
with group law given by the star-product~\eqref{int.1}. Then
$G^{-\infty}_{\susn{2n}}(X;E)$ is a retraction of
$G^{-\infty}_{\susn{2n}}(X;E)[[\varepsilon]].$

Our first main result
is the following.

\begin{theorem_int}\label{Theorem1}
There is a multiplicative `adiabatic' determinant function
\begin{equation*}
\begin{gathered}
\deta: G^{-\infty}_{\susn{2n}}(X;E)[[\varepsilon]]\longrightarrow \bbC^{*},  \\
\deta(A\fsp B)= \deta(A)\deta(B)\
\forall\ A,B\in G^{-\infty}_{\susn{2n}}(X;E)[[\varepsilon]],   
\end{gathered}
\end{equation*}
which generates $\coH^1(G^{-\infty}_{\susn{2n}}(X;E)[[\varepsilon]]).$
\end{theorem_int}
\noindent This is proven in \S\ref{adiabatic.0} by considering a corresponding
determinant for mixed isotropic operators and taking the adiabatic limit.

Given a (locally trivial) fibration of compact manifolds
\begin{equation}
\xymatrix{Z\ar@{-}[r]&M\ar[d]^{\phi}\\&B}
\label{int.4}\end{equation}
and a family of elliptic $2n$-suspended operators 
\begin{equation*}
              D\in \Psi_{\susn{2n}}^{k}(M/B;E,F)
\label{perdet.137}\end{equation*}
with vanishing numerical index, one can construct an associated
determinant line bundle $\Deta(D)\to B$ as described in 
\S\ref{adiabatic.0}, the definition being in terms of (a slightly extended
notion of) principal bundles; a related construction can be found in
\cite{Pressley-Segal1}. More generally, this construction can be extended to
a fully elliptic family of product-suspended operators (see the appendix
and \S\ref{Psop.0} for the definition) 
\begin{equation*}
              D\in \Psi_{\psusn{2n}}^{k,k'}(M/B;E,F)
\end{equation*}
with vanishing numerical index. Our second result is to relate this 
determinant line bundle with Quillen's definition via Bott periodicity.
Let $D_0\in \Psi^{1}(M/B;E,F)$ be a family of elliptic operators with 
vanishing numerical index. Define, by recurrence for $n\in\bbN,$ the
fully elliptic product-suspended familes
by
\begin{equation*}
D_{n}(t_{1}, \ldots, t_{n}, \tau_{1},\ldots, \tau_{n})=
\begin{pmatrix}
              it_{n}-\tau_{n}  & D_{n-1}^{*}  \\
                  D_{n-1}   & it_{n}+\tau_{n}
\end{pmatrix}\in\Psi_{\psusn{2n}}^{1,1}(M/B;2^{n-1}(E\oplus F)),
\label{int.5}\end{equation*}
where $2^{n-1}(E\oplus F)$ is the direct sum of $2^{n-1}$ copies of
$E\oplus F.$ In \S\ref{detn.0} we prove

\begin{theorem_int}[Periodicity of the determinant line bundle]\label{int.6}
For each $n\in \bbN$, their is an isomorphism
$\Deta(D_{n})\cong\Det(D_0)$ as line bundles over $B.$
\end{theorem_int}

In \S\ref{eta.0}, we investigate the counterpart of the eta invariant for
the determinant of theorem~\ref{Theorem1}. After extending the definition
given in \cite{MR96h:58169} to product-suspended operators, we relate this
invariant (denoted here $\eta_{\sus})$ to the extension of the original
spectral definition of Atiyah, Patodi and Singer given by Wodzicki
\cite{Wodzicki3}. Namely consider
\begin{equation}
\eta_z(A)=\sum_{j} \sgn (a_{j}) | a_{j} |^{-z}
\label{perdet.108}\end{equation}
where the $a_{j}$ are the eigenvalues of $A$ in order of increasing $|a_{j}|$
repeated with multiplicity.

\begin{theorem_int} If $A\in\Psi^{1}(X;E)$ is an invertible self-adjoint elliptic
pseudodifferential operator and $A(\tau)=A+i\tau\in\Psi^{1,1}_{\psus}(X;E)$
is the corresponding product-suspended family then  
\begin{equation}
\eta_{\sus}(A(\tau))=\reg_{z=0}\eta_z(A)=\eta (A)
\label{perdet.109}\end{equation}
is the regularized value at $z=0$ of the analytic extension of
\eqref{perdet.108} from its domain of convergence.
\label{int.7}\end{theorem_int}

The eta invariant for product-suspended operators is, as in the suspended
case discussed in (\cite{MR96h:58169}), a log-multiplicative functional  
\begin{equation*}
          \eta_{\sus}(AB)= \eta_{\sus}(A)+\eta_{\sus}(B),\
   A\in \Psi_{\psus}^{k,k'}(X;E), \; B \in \Psi_{\psus}^{l,l'}(X;E).
\label{perdet.110}\end{equation*}

Finally, in \S\ref{etat.0}, we show (see Theorem~\ref{perdet.81}) that in the
appropiate context, this eta invariant can be interpreted as the logarithm
of the determinant of Theorem~\ref{Theorem1}.  

To discuss these results, substantial use is made of various classes of
pseudodifferential operators, in particular product-type suspended
operators and mixed isotropic operators.  An overview of the various
classes used in this paper is given in \S\ref{Psop.0} and some of their
properties are discussed in the appendix.

\paperbody

\section{Determinant line bundle}
\label{detdirac.0}

Quillen in \cite{Quillen} introduced the determinant line bundle for a
family of $\overline{\pa}$ operators. Shortly after, Bismut and Freed in
\cite{Bismut-Freed} and \cite{Bismut-Freed1} generalized the definition to
Dirac operators.  We will show here that this is induced by the Fredholm
determinant, as a representation of the group $G^{-\infty}.$ To do so we
need to slightly generalize the standard notion of a principal bundle.

\subsection{Bundles of groups}

\begin{definition}\label{gpb.1}
Let $G$ be a topological group (possibly infinite dimensional). Then a
fibration $\mathcal{G}\to B$ over a compact manifold $B$ with typical fibre
$G$ is called a \textbf{bundle of groups with model $G$} if its structure
group is contained in $\Aut(G),$ the group of automorphisms of $G.$
\end{definition}

The main example of interest here is the bundle of smoothing groups, with fibre
$G^{-\infty}(Z_b)$ on the fibres of a fibration \eqref{int.4}. In this
case the group is smooth and the bundle inherits a smooth structure.

\begin{definition}\label{gpb.2} Let $\phi:\cG\longrightarrow B$ be a bundle
of groups with model $G,$ then a (right) \textbf{principal $\cG$-bundle}  
is a smooth fibration $\pi:\cP\longrightarrow B$ with typical fibre $G$ 
together with a continuous (or smooth) fibrewise group action 
\begin{equation*}
h:\cP_b\times\cG_b\ni(p,g)\longmapsto p\cdot g^{-1}\in\cP_b 
\label{perdet.91}\end{equation*}
which is continous (or smooth) in all variables, locally trivial and free and
transitive on the fibres. An isomorphism of principal $\cG$-bundles is an
isomorphism of the total spaces which intertwines the group actions.
\end{definition}

The fibre actions combine to give a continous map from the fibre product
\begin{equation*}
\cP\times_{B}\cG=\{(p,g)\in \mathcal{P}\times \mathcal{G};
\pi(p)= \phi(g)\}\longrightarrow \cP.
\label{perdet.92}\end{equation*}

Definition~\ref{gpb.2} is a generalization of the usual notion of a
principal bundle for a group $G$ in the sense that a principal $G$-bundle
$\pi:\mathcal{P}\longrightarrow B$ is naturally a principal $\cG$-bundle
for the trivial bundle of groups $\cG=G\times B\to B.$ Any bundle of groups
$\mathcal{G}\to B$ is itself a principal $\mathcal{G}$-bundle and should be
thought of as the trivial principal $\mathcal{G}$-bundle. Thus a principal
$\mathcal{G}$-bundle $\mathcal{P}\to B$ is trivial, as a principal
$\cG$-bundle, if it is isomorphic as a principal $\cG$-bundle to $\cG.$

\subsection{Classifying principal bundles}

\begin{lemma}\label{perdet.94} If $G$ has a topological classifying
sequence \emph{of groups}
\begin{equation}
G\longrightarrow\EG\longrightarrow\BG
\label{perdet.95}\end{equation}
(so $\EG$ is weakly contractible) which is a Serre fibration, $\cG$ is a
bundle of groups modelled on $G$ with structure group $H\subset\Aut(EG,G),$
the group of automorphisms of $EG$ restricting to automorphisms of $G,$ 
then, principal $\cG$-bundles over compact bases are classified up to
$\cG$-isomorphism by homotoply classes of global sections of a bundle
$\cG(\BG)$ of groups with typical fibre $\BG.$
\end{lemma}

\begin{proof} The assumption that the structure group of $\cG$ is a
subgroup of $\Aut(EG,G)$ allows the bundle of groups $\cG$ to be extended
to a bundle of groups with model $\EG.$ Namely taking an open cover of $X$
by sets over which $\cG$ is trivial, the fibres may be extended to $\EG,$
the transition maps then extend to the larger fibres and the cocycle condition
continues to hold. Denote the resulting bundle of groups,
$\cG(\EG)\supset\cG,$ with typical fibre $\EG.$ The quotiet bundle 
\begin{equation*}
\cG(\BG)=\cG(\EG)/\cG
\label{perdet.113}\end{equation*}
is a bundle of groups with typical fibre $\BG$ and structure group
$\Aut(\EG,G)$ acting on $\BG.$

Similarly, any (right) principal $\cG$-bundle, $\cP,$ has an extension to a
principal $\cG(\EG)$-bundle, $\cP(\EG),$ 
\begin{equation*}
\cP(\EG)_x=\cP_x\times\cG(\EG)_x/\cG_x,\ (p,e)\equiv (pg^{-1},eg^{-1}).
\label{perdet.99}\end{equation*}
Since the group $\EG$ is, by hypothesis, weakly contractible, and the base
is compact, the extended bundle $\cP(\EG)$ has a continuous global
section. As in the case of a traditional principal bundle, the quotient of
this section by the fibrewise action of $\cG$ gives a section of
$\cG(\BG).$ Since 
all sections of a bundle with contractible fibre are homotopic, the section
of $\cG(\BG)$ is well-defined up to homotopy. Bundles isomorphic as principal
$\cG$ bundles give homotopic sections and the construction can be reversed
as in the standard case. Namely, given 
a continuous section $u:B\longrightarrow\cG(\BG)$ we may choose a `good' open
cover, $\{U_i\}$ of $B,$ so that each of the open sets is
contractible and $\cG$ is trival over them. By assumption, the sequence
\eqref{perdet.95} is a Serre fibration, and the fibre is weakly
contractible, so it follows that 
$u$ lifts to a global section  $\tilde u:B\longrightarrow\cG(\EG).$ The
subbundle, given by the fibres $G\subset\EG$ in local trivializations, is
well-defined and patches to a principal $\cG$ bundle from which the given
section can be recovered.
\end{proof}

\subsection{Associated bundles}

As in the usual case there is a notion of a vector bundle associated to a
principal $\cG$-bundle. Suppose given a fixed (real or complex) vector
space $V$ and a smooth bundle map  $r: \cG\times V \to B\times V$
which is a family of representations,
\begin{equation*}
     r_{b}:\cG_{b}\times V \to V
\label{perdet.103}\end{equation*}
of the $\cG_{b}.$ Then, from $\cP$ and $r,$ one can form the
\emph{associated vector bundle} $\cP\times_{r} V$
with fibre 
\begin{equation*}
          (\mathcal{P}\times_{r} V)_{b}= \mathcal{P}_{b}\times V / \sim_{b}
\label{perdet.104}\end{equation*}
where $\sim_{b}$ is the equivalence relation 
\begin{equation*}
(pg, r_{b}(g^{-1},v)) \sim_{b} (p,v).
\label{perdet.105}\end{equation*}

\subsection{$\Det(\cP)$}

Consider again the fibration of closed
manifolds \eqref{int.4} and let  
\begin{equation}
D\in\Psi^m(M/B;\bbE),\ D:\CI(M;E^{+})\to \CI(M;E^{-})
\label{fam.1}\end{equation} 
be a family of elliptic operators parametrized by the base $B$. Then
\begin{equation*}
\xymatrix{\cG^{-\infty}(M/B;E^{+})\ar[d]\\ B}
\label{db.1}\end{equation*}
with fibres
\begin{equation*}
G^{-\infty}(Z_{b};E_{b}^{+})=
\left\{\Id +Q;Q\in \Psi^{-\infty}(Z_{b};E^{+}(b)),\ \Id+Q_{b}
 \text{ is invertible}\right\}
\label{perdet.100}\end{equation*}
is a bundle of groups, with model $G^{-\infty}.$ To the family $D$ we
associate the bundle
\begin{equation}
\xymatrix{G^{-\infty}\ar@{-}[r]&\mathcal{P}(D)\ar[d]\\&B}
\label{db.2}\end{equation}
of invertible perturbations of $D$ by smoothing operators where the fibre
at $b$ is
\begin{equation*}
\cP_{b}(D)=\left\{D_{b}+Q_{b};Q_{b}\in\Psi^{-\infty}(Z_{b};E^{+},E^{-}),\
D_{b}+Q_{b}\text{ is invertible}\right\}.
\label{perdet.101}\end{equation*}
The assumption that the numerical index vanishes implies that $\cP_b(D)$ is
non-empty. In fact, for each $b\in B,$ the group
$G^{-\infty}(Z_{b};E^{+}(b))$ acts freely and transitively on the right
on $\mathcal{P}_{b}(D)$ to give $\cP(D)$ the structure of a principal
$\cG^{-\infty}(M/B;E^{+})$-bundle. On 
the other hand, the Fredholm determinant gives a smooth map 
\begin{equation*}
    \det: \mathcal{G}^{-\infty}\longrightarrow \bbC^{*}\cong \GL{1}{\bbC}
\label{perdet.102}\end{equation*}
which restricts to each fibre to a representation.

Thus the construction above gives a line bundle associated to the principal
bundle \eqref{db.2}; for the moment we denote it $\Det(\cP).$

\subsection{Quillen's definition}

\begin{proposition}\label{detdirac.15} For an elliptic family of
pseudodifferential operators with vanishing numerical index, the
determinant line bundle of Quillen, $\Det(D),$ is naturally isomorphic to
the line bundle, $\Det(\cP),$ associated to the bundle \eqref{db.2} and the
determinant as a representation of the structure group.
\end{proposition}

\begin{proof} First we recall Quillen's definition (following Bismut and
Freed \cite{Bismut-Freed1}). Since it extends readily we consider a
pseudodifferential version rather than the original context of Dirac
operators. So, for a fibration as in \eqref{int.4}, let $D$ be the smooth
family of elliptic pseudodifferential operators of \eqref{fam.1}.  Suppose
that it is of positive order.  We also set $E^\pm(b)=E\big|_{Z_b},$
$\cE^{\pm}_b=\CI(Z_b,E^{\pm}(b))$ and consider the infinite dimensional
bundles $\cE^{\pm}$ over $B.$

Assume that $D_{b}$ has vanishing numerical index. Choosing
inner products on $E^{\pm}$ and a positive smooth density on the fibres of
$M$ allows the adjoint $D^*$ of $D$ to be defined. Then, for each $b\in B,$
$D^{*}_{b}D_{b}:\cE_{b}^{+}\longrightarrow \cE_{b}^{+}$ and
$D_{b}D^{*}_{b}:\cE_{b}^{-}\longrightarrow\cE_{b}^{-}$ have a discrete
spectrum with nonnegative eigenvalues. They have the same positive
eigenvalues with $D_b$ an isomorphism of the corresponding
eigenspaces. Given $\lambda>0,$ the sets
\begin{equation*}
\mathcal{U}_{\lambda}=\left\{b\in B;\lambda\text{ is not an eigenvalue of } 
D_{b}^{*}D_{b}\right\} 
\label{detdirac.2}\end{equation*} 
are open and $\cH^{+}_{[0,\lambda)}\subset\cE^{+}$ and
$\cH^{-}_{[0,\lambda)}\subset\cE^{-},$ respectively spanned by the
eigenfunctions of $D^{*}_{b}D_{b}$ and of $D_{b}D^{*}_{b}$ with
eigenvalues less than $\lambda,$ are bundles over $\mathcal{U}_{\lambda}$ of
the same dimension, $k=k(\lambda).$ Now,
$\cH_{[0,\lambda)}=\cH_{[0,\lambda)}^{+}\oplus\cH_{[0,\lambda)}^{-}$ is a 
superbundle to which we associate the local determinant bundle
\begin{equation*}
\Det(\cH_{[0,\lambda)})=(\wedge^{k}\cH_{[0,\lambda)}^{+})^{-1}
\otimes (\wedge^{k}\cH_{[0,\lambda)}^{-}).
\label{detdirac.3}\end{equation*}
A linear map $P:\cH_{[0,\lambda)}^{+}\to \cH_{[0,\lambda)}^{-}$
induces a section 
\begin{equation}
\det(P)=\wedge^{m}P:\wedge^{m} \cH_{[0,\lambda)}^{+}\longrightarrow 
\wedge^{m}\cH_{[0,\lambda)}^{-}
\label{detdirac.4}\end{equation}
of $\Det(\cH_{[0,\lambda)}).$

For $0<\lambda <\mu,$ $\cH_{[0,\mu)}= \cH_{[0,\lambda)}\oplus\cH_{(\lambda,\mu)}$
over $\mathcal{U}_{\lambda}\cap\mathcal{U}_{\mu},$ where 
$\cH_{(\lambda,\mu)}=\cH^{+}_{(\lambda,\mu)}\oplus 
\cH^{-}_{(\lambda,\mu)} $ and $\cH^{+}_{(\lambda,\mu)}$ 
and $\cH^{-}_{(\lambda,\mu)}$ are respectively the local vector 
bundles spanned by the eigenfunctions of $D^{*}_{b}D_{b}$ and
$D_{b}D^{*}_{b}$ with associated eigenvalues between
$\lambda$ and $\mu.$ Thus, if $D_{(\lambda,\mu)}$ denotes the restriction of 
$D$ to $\cH^{+}_{(\lambda,\mu)}$, then \eqref{detdirac.4} leads to
transition maps
\begin{equation*}
\phi_{\lambda,\mu}:\Det(\cH_{[0,\lambda)})\ni s\longmapsto
  s\otimes\det(D_{(\lambda,\mu)})\in\Det(\cH_{[0,\mu)})\text{ over }
\mathcal{U}_{\lambda}\cap\mathcal{U}_{\mu}.
\label{detdirac.7}\end{equation*}
The cocycle conditions holds over triple intersections and the resulting
bundle, which is independent of choices made (up to natural isomorphism),
is Quillen's determinant bundle, $\Det(D).$ 

Let $Q_{b}\in \Psi^{-\infty}(Z_{b};E^{+},E^{-}),$ for
$b\in \mathcal{U}\subset B$ open, be a smooth family of perturbations such that 
$D_{b}+Q_{b}$ is invertible; it therefore gives a section of $\mathcal{P}$
over $\cU.$ The associated bundle $\Det(\cP)$ is then also trivial over
$\cU$ with
\begin{equation*}
\cU\ni b\longrightarrow (D_b+Q_b,1)
\label{perdet.114}\end{equation*}
being a non-vanishing section. For $\lambda>0$, let $P_{[0,\lambda)}$
be the projection  
onto $\cH_{[0,\lambda)}$, and denote by $P_{[0,\lambda)}^{+}$ and 
$P_{[0,\lambda)}^{-}$ the projections onto $\cH_{[0,\lambda)}^{+}$
and $\cH_{[0,\lambda)}^{-}$ respectively.  Then, on $\mathcal{U}\cap
\mathcal{U}_{\lambda}$ for $\lambda$ large enough,
$P^{-}_{[0,\lambda)}(D_{b}+Q_{b})P^{+}_{[0,\lambda)}$ is invertible, and one
can associate to the section $D_{b}+Q_{b}$ of $\mathcal{P}$ the isomorphism
\begin{multline}
F_{\mathcal{U},\lambda}:\Det(\mathcal{P})\ni
[(D_{b}+Q_{b},c)]\longmapsto\\
\det(P^{-}_{[0,\lambda)}
(D_{b}+Q_{b})P^{+}_{[0,\lambda)})\det(A(Q_{b},\lambda))c\in\Det(D),
\label{detdirac.12}\end{multline}
where $\det(P^{-}_{[0,\lambda)}(D_{b}+Q_{b})P^{+}_{[0,\lambda)})$ is defined 
by \eqref{detdirac.4},
\begin{equation}
A(Q_{b},\lambda)= (D_{b}+P^{-}_{[0,\lambda)}Q_{b}P^{+}_{[0,\lambda)})^{-1}
(D_{b} +Q_{b})\in G_{b}^{-\infty},
\label{detdirac.13}\end{equation}
and $\det(A(Q_{b},\lambda))\in\bbC^{*}$ is the determinant defined on 
$G_{b}^{-\infty}.$

The map $\mathcal{F}_{\mathcal{U},\lambda}$ induces a global isomorphism of
the two notions of determinant bundle since it is independent of
choices. Indeed, it is compatible with the equivalence relation $\sim_{b}$ in the
sense for each $g\in G^{-\infty}(Z_{b};\cE^{+})$ such that both
$P^{-}_{[0,\lambda)}(D_{b}+Q_{b})P^{+}_{[0,\lambda)}$ and
    $P^{-}_{[0,\lambda)}(D_{b}+Q_{b})gP^{+}_{[0,\lambda)}$ are invertible,
\begin{equation*}
     \mathcal{F}_{\mathcal{U},\lambda}( (D+Q_{b})g,\det(g^{-1})c)= 
     \mathcal{F}_{\mathcal{U},\lambda}( (D+Q_{b}), c).
\label{perdet.107}\end{equation*}
It is also compatible with increase of $\lambda$ to $\mu$ in that
$\phi_{\lambda,\mu}\circ F_{\mathcal{U},\lambda}= F_{\mathcal{U},\mu}$ on
$\mathcal{U}\cap\mathcal{U}_{\lambda}\cap\mathcal{U}_{\mu}.$ This is
readily checked
\begin{equation}
\begin{aligned}
\phi_{\lambda,\mu}\circ F_{\mathcal{U},\lambda}(D_{b}+Q_{b},c)&=
\phi_{\lambda,\mu}[\det(P^{-}_{[0,\lambda)}(D_{b}+
Q_{b})P^{+}_{[0,\lambda)}) \det(A(Q_{b},
\lambda))c]  \\
&=\det(A(Q_{b},\lambda))\det(P^{-}_{[0,\lambda)}(D_{b}+
Q_{b})P^{+}_{[0,\lambda)})\otimes\det(D^{+}_{(\lambda,\mu)})c \\
&=\det(A(Q_{b},\lambda))\det(P^{-}_{[0,\mu)}(D_{b}+
P^{-}_{[0,\lambda)}Q_{b}P^{+}_{[0,\lambda)})P^{+}_{[0,\mu)})c \\
&=\det(A(Q_{b},\lambda))\det(P^{-}_{[0,\mu)}(D_{b}+
Q_{b})P^{+}_{[0,\mu)}) \, \, \, \, \times \\
&\det((D_{b}+P^{-}_{[0,\mu)}Q_{b}P^{+}_{[0,\mu)})^{-1}
(D_{b}+P^{-}_{[0,\lambda)}Q_{b}P^{+}_{[0,\lambda)}) ) c \\
&=\det(A(Q_{b},\mu))\det(P^{-}_{[0,\mu)}(D_{b}+
Q_{b})P^{+}_{[0,\mu)}) c\\
&= F_{\mathcal{U},\mu}(D_{b}+Q_{b},c) \, \, . 
\end{aligned}
\label{detdirac.14}\end{equation}
\end{proof}

\subsection{Metric on $\Det(\cP)$}

The Quillen metric has a rather direct expression in terms of the
definition of the determinant bundle as $\Det(\cP).$ Namely, if
$(D_{b}+Q_{b})$ is a section of $\cP$ over the open set $\mathcal{U}\subset
B$, then
\begin{equation}
\left|(D_{b}+Q_{b},1)\right|_{Q}=\exp{(-\frac{1}{2}\zeta'_{b}(0))},
\label{detdirac.16}\end{equation}
where $\zeta_{b}$ is the $\zeta$-function associated to the self-adjoint
positive elliptic operator $(D_{b}+Q_{b})^{*}(D_{b}+Q_{b})$ as constructed
by Seeley \cite{MR38:6220}. When $A_{b}=\Id+R_{b}\in
G_{b}^{-\infty}$  with $R_{b}:\cH_{[0,\lambda)}^{+}\to
\cH_{[0,\lambda)}^{+}$ for some $\lambda>0,$ Proposition
9.36 of \cite{Berline-Getzler-Vergne1}, adapted to this context, shows that 
\begin{equation}
\left|(D_{b}+Q_{b})A_{b}\right|_{Q}=|\det(A_{b})|\left|D_{b}+Q_{b}
\right|_{Q}
\label{detdirac.17}\end{equation}
but then by continuity the same formula follows in general.
Moreover, in the form \eqref{detdirac.16}, Quillen's metric generalizes
immediatly to the case of an arbitrary family of elliptic pseudodiffrential
operators with vanishing numerical index.

\subsection{Primitivity}

\begin{lemma}\label{perdet.115} The deteminant bundle is `primitive' in the
  sense that there is a natural isomorphism 
\begin{equation}
\Det(PQ)\simeq\Det(P)\otimes\Det(Q)
\label{perdet.116}\end{equation}
for any elliptic families $Q\in\Psi^m(M/B;E,F),$ $P\in\Psi^{m'}(M/B;F,G)$
of vanishing numerical index.
\end{lemma}

\begin{proof}
Let $\mathcal{P}$ and $\mathcal{Q}$ denote the principal bundles of invertible
smoothing perturbations of $P$ and $Q.$ Let $P_b+R_b$ and $Q_b+S_b$ be
local smooth sections over some open set $U.$ Certainly
$L_b=(P_b+R_b)(Q_b+S_b)$ is a local section of the principal bundle for $PQ$
and $(L_b,1)$ as a local section of $\Det(PQ)$ may be identified with the
product of the sections $(P_b+R_b,1)$ and $(Q_b,S_b,1)$ as a section of
$\Det(P)\otimes\Det(Q).$ Changing the section of $\cP$ to $(P_b+R_b)g_b$
modifies the section $L_b$ to $L_bg'_b,$ $g'_b=(Q_b+S_b)^{-1}g_b(Q_b+S_b).$
Since  
\begin{equation*}
\det(g'_b)=\det(g_b)
\label{perdet.117}\end{equation*}
the identification is independent of choices of sections and hence is
global and natural.
\end{proof}

Later, it will be convenient to restrict attention to first order elliptic
operators. This is not a strong restriction since for $k\in \bbZ,$ let
$D\in\Psi^{k}(M/B;E,F)$ be a smooth family of elliptic pseudodifferential
operators with vanishing numerical index. Let
$\Delta_{M/B}\in\Psi^{2}(M/B;F)$ be an associated family of Laplacians, so
that $\Delta_{M/B}+\Id$ is a family of invertible operators.

\begin{corollary} 
The family $D'=(\Delta_{M/B}+\Id)^{-\frac{k-1}{2}}D\in\Psi^{1}(M/B;E,F)$ has
determinant bundle isomorphic to the determinant bundle of $D$. 
\label{fior.1}\end{corollary}

\section{Classes of pseudodifferential operators}
\label{Psop.0}
Since several different types, and in particular combinations of types, of
pseudodifferential operators are used here it seems appropriate to quickly
review the essentials.

\subsection{$\Psi^m(X;E,F)$}

On a compact manifold without boundary the `traditional' algebra (so
consisting of `classical' operators) may be defined in two steps using a
quantization map. The smoothing operators acting between two bundles $E$
and $F$ may be identified as the space 
\begin{equation}
\Psi^{-\infty}(X;E,F)=\CI(X^2;\Hom(E,F)\otimes\Omega_R).
\label{fipomb2.6}\end{equation}
Here $\Hom(E,F)_{x,x'}=E_x\otimes F'_{x'}$ is the `big' homomorphism bundle
and $\Omega =\pi_R^*\Omega$ is the lift of the density bundle from the
right factor under the projection $\pi_R:X^2\longrightarrow X.$ The space
$\Psi^m(X;E,F)$ may be identified with the conormal sections, with respect
to the diagonal, of the same bundle 
\begin{equation}
\Psi^m(X;E,F)=I^m_{\cl}(X^2,\Diag;\Hom(E,F)\otimes\Omega_R).
\label{perdet.147}\end{equation}
More explicitly Weyl quantization,
given by the inverse fibre Fourier transform from $T^*X$ to $TX,$
\begin{multline}
q_g:\rho ^{-m}\CI(\com{T^*X};\pi^*\hom(E,F))\ni a\longmapsto\\
(2\pi)^{-n}\int_{T^*X}\chi\exp{(iv(x,y)\cdot\xi)} a(m(x,y),\xi)d\xi dg
\in\Psi^m(X;E,F)
\label{fipomb2.7}\end{multline}
is surjective modulo $\Psi^{-\infty}(X;E,F).$ Here a Riemann metric, $g,$
is chosen on $X$ and used to determine a small geodesically convex
neighbourhood $U$ of the diagonal in $X^2$ which is identified as a
neighbourhood $U'$ of the zero section in $TX$ by mapping $(x,y)\in U$ to
$m(x,y),$ the mid-point of the geodesic joining them in $X$ and to
$v(x,y)\in T_{m(x,y)}X,$ the tangent vector to the geodesic at that
mid-point in terms of the length parameterization of the geodesic from $y$
to $x.$ The cut-off $\chi \in\CIc(U')$ is taken to be identically equal to
$1$ in a smaller neighbourhood of the diagonal. Connections on $E$ and $F$
are chosen and used to identify $\Hom(E,F)$ over $U$ with the lift of
$\hom(E,F)$ to $U',$ $d\xi$ is the fibre density from $g$ on $T^*X$ and
$dg$ is the Riemannian density on the right. The symbol $a$ is a classical
symbol of order $k$ on $T^*X$ realized as $\rho ^{-k}a'$ where
$a'\in\CI(\com{T^*X})$ with $\com{T^*X}$ the compact manifold with boundary
arising from the radial compactification of the fibres of $T^*X$ and $\rho
_g=|\xi|_g^{-1}$ outside a compact set in $T^*X$ is a boundary defining
function for that compactification.

Then $q_g(a)\in\Psi^{-\infty}(X;E,F)$ if and only if $a\in\dCI(\com{T^*X})$
is a smooth function vanishing to all orders on the boundary of
$\com{T^*X},$ \ie is a symbol of order $-\infty.$ This leads to the short
exact `full symbol sequence'
\begin{equation}
\Psi^{-\infty}(X;E,F)\longrightarrow\Psi^\infty(X;E,F)\overset{\sigma_g}
\longrightarrow\CI(S^*X;\hom(E,F))[[\rho ,\rho ^{-1}]]
\label{fipomb2.8}\end{equation}
with values in the Laurent series in $\rho$ (\ie formal power series in
$\rho$ with finite factors of $\rho ^{-1}).$ The leading part of this is
the principal symbol 
\begin{equation}
\Psi^{m-1}(X;E,F)\longrightarrow\Psi^m(X;E,F)\overset{\sigma_m}{\longrightarrow}
\CI(S^*X;\hom(E,F)\otimes R_m)
\label{perdet.159}\end{equation}
where $R_m$ is the trivial bundle with sections which are homogeneous of
degree $m$ over $T^*X\setminus0.$ Pseudodifferential operators act from
$\CI(X;E)$ to $\CI(X;F)$ and composition gives a filtered product,
\begin{equation}
\Psi^m(X;F,G)\circ\Psi^{m'}(X;E,F)\subset\Psi^{m+m'}(X;E,G)
\label{fipomb2.9}\end{equation}
which induces a star product on the image spaces in \eqref{fipomb2.8},  
\begin{equation}
a\star_g b=ab+\sum\limits_{j=1}^\infty B_j(a,b)
\label{fipomb2.10}\end{equation}
where the $B_j$ are smooth bilinear differential operators with polynomial
coefficients on $T^*X$ lowering total order, in terms of power series, by $j.$
The leading part gives the multiplicativity of the principal symbol.

\subsection{$\Psi^{m}_{\susn{p}}(X;E,F)$}

There is a natural Fr\'echet
topology on $\Psi^m(X;E,F),$ corresponding to the \ci\ topology on the
symbol and the kernel away from the diagonal. Thus, smoothness of maps into this
space is well-defined. The $p$-fold suspended operators are a subspace 
\begin{equation}
\Psi^{m}_{\susn{p}}(X;E,F)\subset
\CI\left(\bbR^m;\Psi^{m}(X;E,F)\right)
\label{perdet.148}\end{equation}
in which the parameter-dependence is symbolic (and classical). In terms of
the identification \eqref{perdet.147} this reduces to 
\begin{multline}
\Psi^m_{\susn{p}}(X;E,F)=\\
\cF^{-1}_{\bbR^p}\left(
I^M_{\cl,\cS}(X^2\times\bbR^p,\Diag\times\{0\};\Hom(E,F)\otimes\Omega_R)\right),
\ M=m+\frac p4.
\label{perdet.149}\end{multline}
Here we consider conormal distributions on the non-compact space
$X^2\times\bbR^p$ but with respect to the compact submanifold
$\Diag\times\{0\};$ the suffix $\cS$ denotes that they are to be Schwartz at
infinity and then the inverse Fourier transform is taken in the Euclidean
variables $\bbR^p$ giving the `symbolic' parameters. The shift of $m$ to
$M$ is purely notational. These kernels can also be expressed directly as
in \eqref{fipomb2.7} with $a$ replaced by 
\begin{equation}
a\in\rho ^{-m}\left(\com{T^*X\times\bbR^p};\pi^*\hom(E,F)\right).
\label{perdet.150}\end{equation}
Composition, mapping and symbolic properties are completely analogous to
the `unsuspended' case. Note that we use the abbreviated notation for suffixes
$\susn{1}=\sus.$

If $D$ is a first order elliptic differential operator acting on a bundle
on $X$ then $D+i\tau\in\Psi^1_{\sus}(X;E)$ is elliptic in this sense and
invertible, with inverse in $\Psi^{-1}_{\sus}(X;E),$ if $D$ is self-adjoint
and invertible. However this is not the case for general (elliptic
self-adjoint) $D\in\Psi^1(X;E);$ we therefore introduce larger spaces which
will capture these operators and their inverses.

\subsection{$\Psi^{m,m'}_{\psusn{p}}(X;E,F)$}

By definition in \eqref{perdet.149}, before the inverse Fourier transform
is taken, the singularities of the `kernel' are constrained to
$\Diag\times\{0\}\subset X^2\times\bbR^p.$ For product-type (really
partially-product-type corresponding to the fibration of $X\times\bbR^p$
with base $\bbR^p)$ the singularities are allowed to fill out the larger
submanifold
\begin{equation}
X^2\times\{0\}\supset\Diag\times\{0\}.
\label{perdet.151}\end{equation}
Of course they are not permitted to have arbitrary singularities but rather
to be conormal with respect to these two, nested, submanifolds 
\begin{multline}
\Psi^{m,m'}_{\psusn{p}}(X;E,F)=\\
\cF^{-1}_{\bbR^p}\left(
I^{M'M}_{\cl,\cS}(X^2\times\bbR^p,X^2\times\{0\},\Diag\times\{0\};
\Hom(E,F)\otimes\Omega_R)\right),\\
\ M=m+\frac p4,\ M'=m'+\frac p4-\frac n2.
\label{perdet.152}\end{multline} 
The space of classical product-type pseudodifferential operators is discussed
succinctly in an appendix below. Away from $\Diag\times\{0\}$ the elements
of the space on the right are just classical conormal distributions at
$\{0\}\times\bbR^p,$ so if $\chi\in\CI(X^2)$ vanishes near the diagonal (or
even just to infinite order on it) 
\begin{equation}
K\in\Psi^{m,m'}_{\psusn{p}}(X;E,F)\Longrightarrow \chi
K\in\rho ^{-m'}\CI(\com{\bbR^p}\times X^2;\Hom(E,F)\otimes\Omega _R)
\label{perdet.153}\end{equation}
is just a classical symbol in the parameters depending smoothly on the
variables in $X^2.$
Conversely, if $\chi'\in\CI(X^2)$ has support sufficiently near the
diagonal then the kernel is given by a formula as in \eqref{fipomb2.7}
\begin{multline}
\chi 'K=(2\pi)^{-n}\int_{T^*X}\chi\exp{(iv(x,y)\cdot\xi)}
a(m(x,y),\xi,\tau)d\xi dg,\\
a\in(\rho'') ^{-m}(\rho')^{-m'}\CI(S;\pi^*\hom(E,F)),\
S=[\com{T^*X\times\bbR^p},0_{T^*X}\times\pa\com{\bbR^p}].
\label{perdet.154}\end{multline}
Here the space on which the `symbols' are smooth functions (apart from the
weight factors) is the same compactification as in \eqref{perdet.150} but
then blown up at the part of the boundary (\ie infinity) corresponding to
finite points in the cotangent bundle. Then $\rho''$ is a defining function
for the `old' part of the boundary and $\rho '$ for the new part, produced
by the blow-up. Conversely \eqref{perdet.154} and \eqref{perdet.153}
together (for a partition of unity) define the space of kernels.

From the general properties of blow-up, if
$\rho\in\CI(\com{T^*X\times\bbR^p})$ is a defining function for the
boundary then $\rho =\rho '\rho ''$ after blow-up. From this it follows
easily that 
\begin{equation}
\Psi^m_{\susn{p}}(X;E,F)\subset \Psi^{m,m}_{\psusn{p}}(X;E,F).
\label{perdet.155}\end{equation}

Again these `product suspended' operators act from $\cS(X\times\bbR^p;E)$
to $\cS(X\times\bbR^p;F)$ and have a doubly-filtered composition 
\begin{equation}
\Psi^{m_1,m_1'}_{\psusn{p}}(X;F,G)\circ\Psi^{m_2,m_2'}_{\psusn{p}}(X;E,F)\subset
\Psi^{m_1+m_2,m_1'+m_2'}_{\psusn{p}}(X;E,G).
\label{perdet.156}\end{equation}
The symbol map remains, but now only corresponds to the part of the
amplitude in \eqref{perdet.154} at $\rho ''=0$  
\begin{equation}
\Psi^{m-1,m'}_{\psusn{p}}(X;E,F)\longrightarrow\Psi^{m,m'}_{\psusn{p}}(X;E,F)
\overset{\sigma_m}\longrightarrow \cS^{m,m'}_{\psusn{p}}(X;E,F)
\label{perdet.157}\end{equation}
with 
\begin{equation*}
  \cS^{m,m'}_{\psusn{d}}(X;E,F)=\CI([S(T^*\times\bbR^p),0\times\bbS^{p-1}];\hom(E,F)\otimes R_{m,m'})
\end{equation*}
the space of smooth sections of a bundle over the sphere bundle
corresponding to $T^*X\times\bbR^p,$ blown up at the image of the zero
section and with $R_{m,m'}$ a trivial bundle capturing the weight factors.

The other part of the amplitude corresponds to a more global `symbol map'
called here the `base family'
\begin{equation}
\Psi^{m,m'-1}_{\psusn{p}}(X;E,F)\longrightarrow\Psi^{m,m'}_{\psusn{p}}(X;E,F)
 \overset{\beta_{m'}}\longrightarrow
\CI(\bbS^{p-1};\Psi^{m}(X;E,F)\otimes R_{m'})
\label{perdet.158}\end{equation}
taking values in pseudodifferential operators on $X$ depending smoothly on
the parameters `at infinity', \ie in $\bbS^{p-1}$ with the appropriate
homogeneity bundle (over $\bbS^{p-1}).$

These two symbol maps are separately surjective and jointly surjective onto
pairs satisfying the natural compatibility condition 
\begin{equation}
\sigma _m(\beta _{m'}(A))=\sigma _{m}(A)\big|_{\pa}
\label{perdet.160}\end{equation}
that the symbol family, restricted to the boundary of the space on which it
is defined, is the symbol family of the base family.

An operator in this product-suspended class is `fully elliptic' if both its
symbol and its base family are invertible. If it is also invertible then
its inverse is in the corresponding space with negated orders. An elliptic
suspended operator is automatically fully elliptic when considered as a
product-suspended operators using \eqref{perdet.155}.

\subsection{$\Psi^{m}_{\wq{2n}{\epsilon}}(\bbR^n)$}

The suspension variables for these product-suspended operators are purely
parameters. However, for the adiabatic limit constructions here, on which
the paper relies heavily, we use products which are non-local in the
parameters.

In the trivial case of $X=\{\text{pt}\}$ we are dealing just with symbols
above and the corresponding non-commutative product is the `isotropic'
algebra of operators on symbols on $\bbR^{2n},$ as operators on $\bbR^n,$
for any $n.$ This is variously known as the Weyl algebra or the 
Moyal product (although both often are taken to mean slightly different
things). The isotropic pseudodifferential operators of order $k$
act on the Schwartz
space $\cS(\bbR^n)$ and, using Euclidean Weyl quantization, may be
identified with with the spaces $\rho^{-k}\CI(\com{\bbR^{2n}}).$ Thus, in
terms of their distributional kernels on $\bbR^{2n},$ this space of
operators is given by essentially the same formula as \eqref{fipomb2.7}
\begin{multline}
q_W:\rho ^{-k}\CI(\com{\bbR^{2n}})\ni b\longmapsto\\
q_W(b)(t,t')=(2\pi)^{-n}\int_{\bbR^{n}}
e^{i(t-t')\cdot\tau}b(\frac{t+t'}{2},\tau)d\tau\in\Psi^{-k}_{\iso}(\bbR^n).
\label{fipomb2.11}\end{multline}
This map is discussed extensively in \cite{Hormander5}. In this case $q_W,$
with inverse $\sigma _W,$ is an isomorphism onto 
the algebra and restricts to an isomorphism of the `residual' algebra
$\Psi^{-\infty}_{\iso}(\bbR^n)=q_W(\cS(\bbR^{2n})).$ The corresponding star
product is the Moyal product.

The full product on symbols on $\bbR^{2n}$ may be written explicitly as 
\begin{equation}
a\circ_\omega b(\zeta)=\pi^{-2n}\int_{\bbR^{8n}}
e^{i\xi\cdot\xi '+i\eta \cdot\eta '+ 2i\omega(\xi',\eta')}
a(\zeta +\xi )b(\zeta +\eta )
|\omega _\xi|^n|\omega _\eta |^n 
\label{perdet.161}\end{equation}
where the integrals are not strictly convergent but are well defined as
oscillatory integrals. Here $\omega$ is the standard symplectic form on
$\bbR^{2n}.$ By simply using linear changes of variables, it may be seen
that this product and the more general ones in which $\omega$ is replaced by an
arbitrary non-degenerate antisymmetric bilinear form on $\bbR^{2n}$ are all
isomorphic. In fact the product depends smoothly on $\omega$ as an
antisymmetric bilinear form, even as it becomes degenerate. When
$\omega\equiv0$ the product reduces to the pointwise, commutative, product
of symbols. In fact it is not necessary to assume that the underlying
Euclidean space is even dimensional for this to be true; of course in the
odd-dimensional case the form cannot be non-degenerate and correspondingly
there is always at least one `commutative' variable.

The adiabatic limit here corresponds to replacing the standard symplectic
from $\omega$ by $\epsilon \omega$ and allowing $\epsilon \downarrow0.$ As
already noted, this gives a family of products on the classical symbol
spaces which is smooth in $\epsilon$ and is the commutative product at
$\epsilon =0.$ We denote the resulting smooth family of algebras by
$\Psi^{m}_{\wq{2n}{\epsilon}}(\bbR^n).$

\subsection{$\Psi^{m,m'}_{\wq{2n}{\epsilon}}(X;E,F)$}

Now, we may replace the parameterized product on the product-suspended algebra
by `quantizing it' as in \eqref{perdet.161}, in addition to the composition in
$X$ itself. For the `adiabatic' choice of $\epsilon\omega $ this induces a
one parameter family of quantized products
\begin{equation}
  [0,1]_{\epsilon}\times \Psi^{m_1,m_1'}_{\psusn{2n}}(X;F,G)\times 
              \Psi_{\psusn{2n}}^{m_2,m'_2}(X;E,F) \longrightarrow 
                \Psi_{\psusn{2n}}^{m_1+m_2,m_1'+m_2'}(X;E,G).
\label{ttt.3}\end{equation}
The suspended operators still form a subalgebra. The Taylor series as
$\epsilon \downarrow 0$ given by
\begin{equation}
(A\spe{\epsilon} B)(u) \sim \sum_{k=0}^{\infty} 
 \left.    \frac{(i\epsilon)^{k}}{2^{k} k!} \omega(D_{v}, D_{w}) A(v)B(w)
\right|_{v=w=u}.
\label{ttt.4}\end{equation}

A more complete discussion of product suspended operators and the 
mixed isotropic product may be found in the appendix.

\subsection{$\Psi^{m,m'}_{\psusn{2n}}(X;E,F)[[\epsilon]]$} This is the
space of formal power series in $\epsilon$ with coefficients in
$\Psi^{m,m'}_{\psusn{2n}}(X;E,F).$ The product \eqref{ttt.4} projects to
induce a product 
\begin{equation}
\Psi^{m_1,m'_1}_{\psusn{2n}}(X;F,G)[[\epsilon]]\times
\Psi^{m_2,m'_2}_{\psusn{2n}}(X;E,F)[[\epsilon]]\longrightarrow 
\Psi^{m_1+m_2,m'_1+m'_2}_{\psusn{2n}}(X;E,G)[[\epsilon]]
\label{perdet.164}\end{equation}
which is consistent with the action on formal power series 
\begin{equation*}
\Psi^{m,m'}_{\psusn{2n}}(X;E,F)[[\epsilon]]\ni A:\CI(X;E)[[\epsilon ]]
\longrightarrow \CI(X;F)[[\epsilon ]].
\label{perdet.165}\end{equation*}

\section{Adiabatic determinant}\label{adiabatic.0}

Let $E\longrightarrow X$ a complex vector bundle over a compact manifold
$X.$ Consider the infinite dimensional group 
\begin{equation*}
G^{-\infty}_{\susn{2n}}(X;E)=\{\Id + Q;\ Q\in \Psi_{\susn{2n}}(X;E),\
          \Id+Q\text{ is invertible}\}
\label{perdet.118}\end{equation*}
of invertible 2n-suspended smoothing perturbations of the identity.  A
naive notion of determinant would be given by using the $1$-form 
\begin{equation*}
    d\log d(A) = \Tr_{\susn{2n}} (A^{-1}dA)
\label{perdet.119}\end{equation*}
where 
\begin{equation*}
\Tr_{\susn{2n}}(B)= \frac{1}{(2\pi)^{2n}} \int_{\bbR^{2n}} \Tr(B(t,\tau)) dtd\tau
\label{perdet.120}\end{equation*}
is the regularized trace for suspended operators as defined in
\cite{MR96h:58169}. The putative determinant is then given by
\begin{equation}
d(A)=\exp\left(\int_{0}^{1}\Tr_{\susn{2n}}(\gamma^{-1}\frac{d\gamma}{ds})ds
\right)
\label{detn.1}\end{equation}
where $\gamma: [0,1]\to G^{-\infty}_{\susn{2n}}(X;E)$ is any smooth path
such that $\gamma(0)=\Id$ and $\gamma(1)=A.$ Although $d(A)$ is
multiplicative, it is topologically trivial, in the sense that for any
smooth loop $\gamma:\bbS^{1}\to G^{-\infty}_{\susn{2n}}(X;E),$
one has
\begin{equation}
\int_{\bbS^{1}}\Tr_{\susn{2n}} (\gamma^{-1}\frac{d\gamma}{ds}) ds
= \frac{1}{(2\pi)^{2n}}  \int_{\bbR^{2n}} \left(\int_{\bbS^{1}}
\Tr(\gamma^{-1}\frac{d\gamma}{ds})ds \right)dt\,d\tau
=0.
\end{equation}
So this is not a topological analogue of the usual determinant.

\subsection{Isotropic determinant}

To obtain a determinant which generates the 1-dimensional cohomology, we
instead use the isotropic quantization of \S\ref{iso.0}. At the cost of
slightly deforming the composition law on $G^{-\infty}_{\susn{2n}}(X;E),$
this determinant will be multiplicative as well.

Notice first that because of the canonical identification  
\begin{equation*}
    \Psi_{\psusn{2n}}^{-\infty,-\infty}(X;E)= \Psi_{\susn{2n}}^{-\infty}(X;E) 
\label{perdet.121}\end{equation*}
there is no distinction between $G^{-\infty}(X;E)$ and the group  
\begin{equation*}
G^{-\infty,-\infty}_{\psusn{2n}}(X;E)=
    \{ \Id +Q;\ Q\in \Psi_{\psusn{2n}}^{-\infty,-\infty}(X;E),\ \Id +Q
    \text{ is invertible,} \},
\label{perdet.122}\end{equation*}
so in this context, we can interchangeably think in terms of suspended or 
product-suspended operators. For $\epsilon \in [0,1]$, we use the
$\spe{\epsilon}$-product of Theorem~\ref{fipomb2.16} to define the group
\begin{multline}
G^{-\infty}_{\wq{2n}{\epsilon}}(X;E)=
    \{\Id+Q;\ Q\in \Psi_{\susn{2n}}^{-\infty}(X;E),\ 
        \exists\ P\in \Psi^{0,0}_{\psusn{2n}}(X;E), \\
         P\spe{\epsilon}(\Id+Q)=(\Id+Q)\spe{\epsilon} P= \Id \}.
\label{dete.1}\end{multline}
For $\epsilon=0$, we have the canonical group isomorphism 
\begin{equation*}
      G^{-\infty}_{\wq{2n}{0}}(X;E)= G^{-\infty}_{\susn{2n}}(X;E).
\label{perdet.123}\end{equation*}
On the other hand, for $\epsilon>0$, the group
$G^{-\infty}_{\wq{2n}{\epsilon}}(X;E)$ is isomorphic to $G^{-\infty}$ so
that it is possible to transfer the Fredholm determinant to it.

\begin{proposition}
For $\epsilon>0$, there is a natural multiplicative determinant 
\begin{equation*}
   \det_{\epsilon}(A): G^{-\infty}_{\wq{2n}{\epsilon}}(X;E)\to \bbC^{*}   
\label{perdet.124}\end{equation*}
defined for $A\in G^{-\infty}_{\wq{2n}{\epsilon}}(X;E)$ by 
\begin{equation*}
\det_{\epsilon}(A)=\exp\left(\int_{0}^{1}
\Tr_{\epsilon}(\gamma^{-1} \spe{\epsilon} \frac{d\gamma}{ds} ) ds \right)
\label{perdet.125}\end{equation*}
where $\gamma: [0,1]\to G^{-\infty}_{\wq{2n}{\epsilon}}(X;E)$ is any smooth
path with $\gamma(0)=\Id$ and $\gamma(1)= A$ so
\begin{equation*}
     d\log \det_{\epsilon}(A) = \Tr_{\epsilon}( A^{-1}\spe{\epsilon}dA).
\label{perdet.126}\end{equation*}
\label{dete.2}\end{proposition}

\begin{proof} To show that $\det_{\epsilon}$ is well-defined and
multiplicative, it suffices to show that it reduces to the Fredholm
determinant under a suitable identification of
$G^{-\infty}_{\wq{2n}{\epsilon}}(X;E)$ with $G^{-\infty}.$  From
Appendix~\ref{wq.0} it follows that
$G^{-\infty}_{\wq{2n}{\epsilon}}(X;E)$ acts on $\cS(X\times \bbR^{n};E)$.
Fix a Riemannian metric $g$ on $X$ and a Hermitian metric $h$ on $E.$ Let
$\Delta\in \Psi^{2}(X;E)$ be the corresponding Laplace operator. Then
consider the mixed isotropic operator 
\begin{equation*}
\Box_{\epsilon}=\Delta+\epsilon\sum_{i=1}^{n}\left(-\frac{\pa^{2}}{\pa t_{i}^{2}}
+ t_{i}^{2}  \right) \,\in \Psi_{\wq{2n}{\epsilon}}^{2}(X;E) 
\label{perdet.127}\end{equation*}
where $\sum_{i=1}^{n}\left( -\frac{\pa^{2}}{\pa t_{i}^{2}}
    + t_{i}^{2}  \right)$ is the harmonic oscillator on $\bbR^{n}.$
As an operator acting on $\cS(X\times\bbR^{n};E)$, $\Box_{\epsilon}$ has
a positive discrete spectrum.  Let $\{\lambda_{k}\}_{k\in \bbN}$ be the
eigenvalues, in non-decreasing order, with corresponding eigensections  
\begin{equation*}
\Box_{\epsilon}f_{k}=\lambda_{k} f_{k}, \quad f_{k}\in \cS(X\times\bbR^{n};E)
\label{perdet.128}\end{equation*}
such that $\{f_{k}\}_{k\in \bbN}$ is an orthonormal basis of
$\Ld(X\times\bbR^{n};E).$ This gives an algebra isomorphism 
\begin{equation*}
\mathcal{F}_{\epsilon}:\Psi^{-\infty}_{\wq{2n}{\epsilon}}(X;E)
\ni A\longmapsto \langle f_{i}, A f_{j}\rangle_{\Ld}\in\Psi^{-\infty}.
\label{perdet.129}\end{equation*}
and a corresponding group isomorphism $\mathcal{F}_{\epsilon}:
G^{-\infty}_{\wq{2n}{\epsilon}}(X;E)\to G^{-\infty}.$ Under these
isomorphisms, one has 
\begin{equation*}
     \Tr_{\epsilon}(A)= \Tr(\mathcal{F}_{\epsilon}(A)) 
\label{perdet.130}\end{equation*}
and consequently 
\begin{equation}
    \det_{\epsilon}(\Id+A)= \det_{\Fr}(\mathcal{F}_{\epsilon}(\Id +A)).
\label{iso.24}\end{equation}
\end{proof}

\subsection{Asymptotics of $\det_\epsilon$}

Now, for any $\delta >0$ we can consider the group of sections, 
\begin{multline}
G^{-\infty}_{\iso}([0,\delta]\times\bbR^{2n}\times X;E) =\{A\in\CI([0,\delta
];\Id+\Psi^{-\infty}(\bbR^{2n}\times X));\\
A(\epsilon)\in G^{-\infty}_{\iso,\epsilon }(\bbR^{2n}\times X;E)\ \forall\
\epsilon \in[0,\delta ]\}.
\label{perdet.32}\end{multline}

\begin{proposition}\label{detsing}
The determinant with respect to the $\spe{\epsilon}$ product defines 
\begin{equation}
\tdet:G^{-\infty}_{\iso}([0,\delta]\times\bbR^{2n}\times X;E)\longrightarrow
\CI((0,\delta ])
\label{perdet.33}\end{equation}
which takes the form 
\begin{equation}
\tdet(A)(\epsilon)= \exp{ \left( \sum_{k=0}^{n-1} \epsilon^{k-n}a_{k}(A)\right)}
F_{\epsilon}(A)\ \forall\ A\in
G^{-\infty}_{\iso}([0,\delta]\times\bbR^{2n}\times X;E)
\label{perdet.34}\end{equation}
where 
\begin{equation}
F:G^{-\infty}_{\iso}([0,\delta]\times\bbR^{2n}\times X;E)\ni A\longmapsto
F_{\epsilon}(A)\in \CI([0,\delta ]) 
\label{perdet.35}\end{equation}
 and $a_{k}:G^{-\infty}_{\iso}([0,\delta]\times\bbR^{2n}\times
X;E)\longrightarrow \bbC$ are \ci\ functions and the $a_k$ only depend on the
Taylor series of $A.$
\end{proposition}

\begin{proof} Since the group is open (for each $\epsilon \in[0,1]$ and
  also for the whole group) the tangent space at any point is simply
  $\CI([0,\delta ];\Psi^{-\infty}(\bbR^{2n}\times X;E)).$ With the usual
  identifications for a Lie group the form  $A^{-1}\spe{\epsilon}dA$
  therefore takes values in $\CI([0,\delta ];\Psi^{-\infty}(\bbR^{2n}\times
  X;E)).$ On the other hand, the trace functional is \emph{not} smooth down
  to $\epsilon =0.$ In fact it is rescaled by a factor of $\epsilon
  ^{-n}.$ Thus,
\begin{equation}
d\log{{\det}_{\epsilon}(A)}=\Tr_{\epsilon}(A^{-1}\spe{\epsilon}dA)\sim
\sum_{k=0}^{\infty} \alpha_{k}\epsilon^{k-n}
\label{iso.26}\end{equation}
is $\epsilon ^{-n}$ times a smooth function. 

For any smooth map  
\begin{equation*}
f: \bbS^{1}\to 
G^{-\infty}_{\iso}([0,\delta]\times\bbR^{2n}\times X;E),
\label{perdet.132}\end{equation*}
the integral $\int_{\bbS^{1}}f^{*}d\log{{\det}_{\epsilon}(A)}$ also has an
asymptotic expansion
\begin{equation}
\int_{\bbS^{1}}f^{*}d\log{{\det}_{\epsilon}(A)} \sim
\sum_{k=0}^{\infty} c_{k}\epsilon^{k-n}\, \, \, , 
c_{k}=\int_{\bbS^{1}}\alpha_{k}\in \bbC.
\label{iso.27}\end{equation} 
On the other hand, by \eqref{iso.24}, this is a winding number so
cannot depend on $\epsilon.$  Hence
\begin{equation}
c_{k}=\int_{\bbS^{1}}\alpha_{k}= 0\Mfor k\ne n.
\label{iso.28}\end{equation}
So, for $k\ne n,$ $\alpha_{k}$ is exact and then \eqref{perdet.34} follows
directly by integration along any path $\gamma:[0,1]\to
G^{-\infty}_{\iso}([0,\delta]\times\bbR^{2n}\times X;E)$ with
$\gamma(0)=\Id$ 
and $\gamma(1)=A.$ The range space is path-connected, so
\begin{equation*}
a_{k}(A)= \int_{0}^{1}\gamma^{*}\alpha_{k}, \quad k<n
\label{perdet.36}\end{equation*}
is independent of the path and well-defined.
\end{proof}

\subsection{Star product}

The restriction map at $\epsilon =0$  
\begin{equation}
R:G^{-\infty}_{\iso}([0,\delta]\times\bbR^{2n}\times X;E)\longrightarrow
G^{-\infty}_{\susn{2n}}(X;E)
\label{perdet.37}\end{equation}
is surjective. From this it follows that if we let $\dot
G^{-\infty}_{\iso}([0,\delta]\times\bbR^{2n}\times X;E)$ be the subgroup of
those elements which are equal to the identity to infinite order at
$\epsilon =0$ then the quotient
\begin{multline}
G^{-\infty}_{\iso}([0,\delta]\times\bbR^{2n}\times X;E)/
\dot G^{-\infty}_{\iso}([0,\delta]\times\bbR^{2n}\times X;E)=\\
G^{-\infty}_{\susn{2n}}(X;E)+\
\epsilon \Psi^{-\infty}(\bbR^{2n}\times X;E)[[\epsilon ]]
\label{perdet.38}\end{multline}
is the obvious formal power series group, namely with invertible leading
term and arbitrary smoothing lower order terms. The composition law is the one
induced by the $\spe{\epsilon}$-product.  Since the higher order terms in
$\epsilon$  amount to an affine extension of the leading part, this
formal power series group is also a classifying group for odd K-theory.

\begin{definition}\label{sus.17b}
We denote by $\Psi^{k,k'}_{\psus(2n)}(X;E)[[\varepsilon]]$, 
$k,k'\in\bbR\cup\{-\infty\}$,
the space of formal series  
\begin{equation*}
    \sum_{\mu=0}^{\infty}a_{\mu}\varepsilon^{\mu}
\label{perdet.133}\end{equation*}
with coefficients $a_{\mu}\in \Psi_{\psusn{2n}}^{k,k'}(X;E)$, 
where $\varepsilon$ is
a formal parameter.
\end{definition}

For $A\in \Psi^{k,k'}_{\psusn{2n}}(X;F,G)[[\varepsilon]]$ 
and $B\in \Psi^{l,l'}_{\psusn{2n}}(X;E,F)[[\varepsilon]]$ the $\fsp$-product
$A \fsp  B\in\Psi^{k+l,k'+l'}_{\psusn{2n}}(X;E,G)[[\varepsilon]]$ is
\begin{equation*}
\begin{aligned}
A\fsp B(u) &= (\sum_{\mu=0}^{\infty}a_{\mu}\varepsilon^{\mu}) \fsp 
          (\sum_{\nu=0}^{\infty}b_{\nu}\varepsilon^{\nu}) \\
       &=\sum_{\mu=0}^{\infty}\sum_{\nu=0}^{\infty}\varepsilon^{\mu+\nu}
           \left( \sum_{p=0}^{\infty}\left.\frac{i^{p}\varepsilon^{p}}{2^{p}p!}
              \omega(D_{v},D_{w})^{p}A(v)B(w)\right|_{v=w=u}\right)
\end{aligned}
\label{perdet.134}\end{equation*}
where $u,v,w\in \bbR^{2n}.$

Since this bases on the asymptotic expansion \eqref{iso.17} of
Appendix~\ref{iso.0}, its associativity follows immediatly from the
associativity of the $\spe{\epsilon}$-product.

\subsection{Adiabatic determinant}

This product is consistent with that of the quotient group in 
\eqref{perdet.38}, so

\begin{lemma}
The quotient group $G^{-\infty}_{\susn{2n}}(X;E)[[\epsilon]]$ of 
\eqref{perdet.38} is canonically isomorphic to
\begin{multline*}
G^{-\infty}_{\susn{2n}}(X;E)[[\varepsilon]]=
\{ (\Id+Q); Q\in \Psi_{\susn{2n}}^{-\infty}(X;E)[[\varepsilon]],\
\exists\ P\in \Psi_{\susn{2n}}^{-\infty}(X;E)[[\varepsilon]]\\ 
\text{ such that }
(\Id+Q)\fsp(\Id+P)=\Id\in\Psi_{\susn{2n}}^{0}(X;E)[[\varepsilon]] \}.
\end{multline*}
\label{sus.19}\end{lemma}

We can now prove Theorem~\ref{Theorem1} stated in the introduction.

\begin{theorem}\label{perdet.39} The functional  
\begin{equation*}
F_{0}(A):G^{-\infty}_{\iso}([0,\delta]\times\bbR^{2n}\times X;E)\to \bbC^*
\label{perdet.135}\end{equation*}
induces a multiplicative determinant $\adet$ on the formal power series group  
\begin{equation*}
G^{-\infty}_{\susn{2n}}(X;E)[[\varepsilon]]
\label{perdet.136}\end{equation*}
in the sense discussed above, \ie it is a smooth multiplicative
function which generates $\coH^1.$
\end{theorem}

\begin{proof}
From \eqref{iso.26},
\begin{equation}
         F_{0}(A)=\exp\left( \int_{\gamma}\alpha_{n}   \right),
\label{detnp.1}\end{equation}
where $\gamma:[0,1]\to G^{-\infty}_{\susn{2n}}(X;E)[[\varepsilon]]$ is any 
smooth path with $\gamma(0)=\Id$ and $\gamma(1)=A$.
In the expansion~\eqref{iso.26}, the only non-trivial cohomological 
contribution comes from $\alpha_{n}.$ Since $\det_{\epsilon}$ 
corresponds to the Fredholm determinant under the identification
of $G^{-\infty}_{\wq{2n}{\epsilon}}(X;E)$ with $G^{-\infty}$ the integral
of $\alpha_{n}$ along a generator of the fundamental group is $\pm 2\pi i.$ Thus,
the determinant induced by $F_{0}(A)$ has the desired topological behavior.  

For the multiplicativity, from \eqref{iso.26},
\begin{equation}
\begin{aligned}
\Tr((A B)^{-1}d(A B))&= \Tr(B^{-1} A^{-1} dA B + B^{-1} A^{-1} A
 dB) \\
   &= \Tr(A^{-1} dA) + 
         \Tr(B^{-1} dB) \, \,,
\end{aligned}
\label{sus.26}\end{equation}
where the $\fsp $-product is used to compose elements and define the inverses.
From the $\varepsilon$-expansion of \eqref{sus.26},
\begin{equation}
\alpha_{n}(A\fsp B)=\alpha_{n}(A) +\alpha_{n}(B).
\label{sus.27}\end{equation}
As a consequence, the determinant defined in \eqref{detn.1} is multiplicative.
\end{proof}

This determinant can be used to define the determinant line bundle 
of a fully elliptic family $D\in
\Psi^{k,k}_{\psus(2n)}(M/B;E,F)[[\varepsilon]]$ of 
fibrewise product $2n$-suspended pseudodifferential operators on a
fibration \eqref{int.4}. Full ellipticity here corresponds to ellipticity of the
leading term $D_{0}$ and its invertibility for large values of the
parameters. Assume in addition that for each $b\in B,$
$D_{b}\in\Psi_{\psusn{2n}}^{k,k}(Z_{b};E_{b},F_{b})[[\varepsilon]]$ can be
perturbed by $Q_{b}\in 
\Psi^{-\infty}_{\susn{2n}}(Z_{b};E_{b},F_{b})[[\varepsilon]]$ to be
invertible, where invertibility is equivalent to invertibility of the
leading term. Then over the manifold $B,$ consider the bundle of invertible
smoothing perturbations with fibres
\begin{multline}
\mathcal{P}_{b}(D)=\{D_{b}+ Q_{b}; Q_{b}\in 
\Psi_{\susn{2n}}^{-\infty}(Z_{b};E_{b},F_{b})[[\varepsilon]], \\
\exists\ P\in\Psi_{\psusn{2n}}^{-k,-k}(Z_{b};F_{b},E_{b})[[\varepsilon]],\
P\fsp (D_{b}+Q_{b})= (D_{b}+Q_{b}) \fsp P =\Id\}.
\label{sus.6}\end{multline}
Let 
\begin{equation}
\xymatrix{ G^{-\infty}_{\susn{2n}}\ar@{-}\ar[r]&
G^{-\infty}_{\susn{2n}}(M/B;E)\ar[d]^{\phi}\\&B} 
\label{bg.15}\end{equation}
be the bundle of groups with fibre at $b\in B$
\begin{multline}  
G^{-\infty}_{\susn{2n}}(Z_{b};E_{b})[[\varepsilon]]= \{ \Id + Q;
Q\in \Psi_{\susn{2n}}^{-\infty}(Z_{b};E_{b})[[\varepsilon]],\\
\exists\ P\in\Psi_{\psusn{2n}}^{0,0}(Z_{b};E_{b})[[\varepsilon]]  
      P\fsp (\Id+Q) = (\Id+Q)\fsp P =\Id \}.
\label{bg.16}\end{multline}
Then $\mathcal{P}(D)$ is a principal
$G^{-\infty}_{\susn{2n}}(M/B;E)[[\varepsilon]]$-bundle in the sense of
Definition~\ref{gpb.2}.

\begin{definition}\label{sus.7}
The \emph{adiabatic determinant line bundle} associated to the family $D$
of product $2n$-suspended elliptic pseudodifferential operators is
\begin{equation*}
      \Deta(D)=\mathcal{P}(D)
      \times_{\det}\bbC
\end{equation*}
induced by the adiabatic determinant as representation of the bundle of
groups \eqref{bg.15}.
\end{definition} 

\section{Periodicity of the numerical index}

In the next section, we establish a relation between the determinant line
bundles of a family of standard elliptic pseudodifferential operators with
the determinant line bundle just defined for families of $2n$-suspended
operators. Here we consider the corresponding question for the numerical index.

\subsection{Product-suspended index}

A product-suspended operator 
\[
P\in\Psi^{m,m'}_{\psusn{k}}(Z;E,F)
\] 
is \textbf{fully elliptic} if both its symbol in the usual sense and its
base family are invertible. Here the base family, elliptic because of the
invertibility of the symbol, is parameterized by $\bbS^{k-1}.$ As a family
of operators over $\bbR^k,$ $P$ has a families index. Since by assumption
the family is invertible at, and hence near, infinity the family defines an
index class in compactly-supported K-theory
\begin{equation}
\ind(P)\in K^0(\bbR^k)=\begin{cases}\bbZ& k\ \even\\ \{0\}& k\ \odd.
\end{cases}
\label{perdet.138}\end{equation}
Thus by choosing a generator (\ie Bott element) in $K^0(\bbR^{2n})$ a
product $2n$-suspended family has a numerical index which we will denote
$\ind_{\susn{n}}$ (since it only arises for even numbers of
parameters). The families index of Atiyah and Singer does not apply
directly to this setting although it \emph{does} apply if the operator is
in the `suspended' subspace (and so in particular $m'=m.)$ Using the
properties of the suspended eta invariant we will show in \S\ref{Adeta}
that the suspended index can be expressed in terms of the `adiabatic' $\eta$ 
invariant discussed below. Namely, suppose a linear decomposition
$\bbR^{2n}=\bbR\times\bbR^{2n-1}$ is chosen in which the variables are
written $\tau$ and $\xi.$ Then, for some $R\in\bbR,$ $P(\tau,\xi)$ is
invertible for $|\tau|\ge R$ for all $\xi\in\bbR^{2n-1}.$ Furthermore, by
standard index arguments we may find a family of smoothing operators, $A,$ of
compact support in $(\tau,\xi)$ such that $P'=P+A$ is invertible for all
$\tau\le R.$ Then
\begin{equation}
\ind_{\susn{n}}(P)=-\frac12\left(\eta_{\adn{n-1}}(P\big|_{\tau=R})
-\eta_{\adn{n-1}}(P'\big|_{\tau=R})\right).
\label{perdet.141}\end{equation}

\subsection{Periodicity}

Here we show that there is a `Bott map' from ordinary pseudodifferential
operators into product-type suspended operators which maps the usual index to the
suspended index (although most of the proof is postponed until later). Thus if
$D\in\Psi^{1}(Z;E,F)$ is an elliptic operator then
\begin{equation}
\bbR^{2} \ni (t,\tau)\longmapsto \hat{D}(t,\tau)=
\begin{pmatrix}
it-\tau & D^{*} \\
D & it+\tau  
\end{pmatrix},\ \hat{D}\in \Psi^{1,1}_{\psusn{2}}(Z;E\oplus F)
\label{fior.2}\end{equation}
is an associated twice-suspended fully elliptic operator. 
In \cite{fipomb2}, such a family is realized explicitly as the indicial 
family of a product-suspended cusp operator.
The ellipticity of $\hat{D}$ follows from the fact that 
\begin{equation}
  \hat{D}^{*}\hat{D}=\begin{pmatrix}
D^{*}D+t^{2}+\tau^{2} & 0 \\
0& D D^{*} +t^{2}+\tau^{2} 
\end{pmatrix}\in\Psi_{\psusn{2}}^{2,2}(Z;E\oplus F)
\label{fior.3}\end{equation}
is an elliptic family which is invertible for $t^2+\tau^2>0.$

\begin{definition}\label{SuDn}
Given an elliptic operator $D\in\Psi^{1}(Z;E,F),$ we define by recurrence on 
$n\in\bbN_{0},$ elliptic product-suspended operators
$D_{n}\in\Psi_{\psusn{2n}}^{1,1}(Z;2^{n-1}(E\oplus F)$ by
\begin{equation*}
   D_{n}(t_{1},\ldots, t_{n},\tau_{1},\ldots,\tau_{n})=
\begin{pmatrix}
it_{n}-\tau_{n} & D_{n-1}^{*} \\
D_{n-1} & it_{n}+\tau_{n}  
\end{pmatrix}
\end{equation*}  
with $D_{0}=D.$
\label{detnpp.1}\end{definition}

\begin{lemma}\label{detn.2}
If $D$ is elliptic then $D_{n}$ is a totally elliptic product $2n$-suspended
operator for all $n$ and $\ind_{\susn{n}}(D_{n})=\ind(D).$
\end{lemma}

\begin{proof} Both the ordinary index and the $n$-suspended index (on
  fully elliptic $2n$-suspended operators) are homotopy invariant. Since
  the map $D\longmapsto D_n$ maps invertible operators to invertible
  operators it follows that $\ind(D)=0$ implies $\ind_{\susn{n}}(D_n)=0.$
  Indeed, $\ind(D)=0$ means there exists a smoothing operator
  $Q\in\Psi^{-\infty}(Z;E,F)$ such that $D+Q$ is invertible. Then
  $(D+sQ)_n$ is a homotopy of fully elliptic $2n$-suspended operators which
  is invertible for $s=1$ so $\ind_{\susn{n}}(D_n)=0.$

The actual equality of the index is proved below in \S\ref{Adeta}, using
\eqref{perdet.141}.
\end{proof}

\section{Periodicity of the determinant line bundle}
\label{detn.0}

\subsection{Adiabatic determinant bundle}

Returning to the setting of a fibration with compact fibres, $\phi:M\to B,$
as in \eqref{int.4}, let $D\in\Psi^{1}(M/B;E,F)$ be a family of elliptic
pseudodifferential operators with vanishing numerical index. From
lemma~\ref{detn.2} (the part that is already proved), the suspended index of the
fully elliptic family $D_n\in\Psi^{1,1}_{\psusn{2n}}(M/B;E,F),$ given by
Definition~\ref{SuDn}, also vanishes. Thus the fibres
\begin{multline}
\mathcal{P}_{\psusn{2n}}(D_{n})_{b}=
\{D_{n,b}+Q_{b}; Q_{b}\in \Psi^{-\infty}_{\susn{2n}}(
Z_{b}; 2^{n-1}(E_{b}\oplus F_{b}))[[\varepsilon]], \\
\exists\ (\hat{D}_{n,b}+Q_{b})^{-1}\in 
\Psi_{\psusn{2n}}^{-1}(Z_{b};2^{n-1}(E_{b}\oplus F_{b}))[[\varepsilon]]\}
 \label{detdirac.31}\end{multline}
are non-empty and combine to give a
principal-$G_{\susn{2n}}^{-\infty}(M/B;2^{n-1}(E\oplus
F))[[\varepsilon]]$-bundle as in \eqref{sus.6}. Since we have defined an
adiabatic determinant on these group we have an associated determinant bundle 
\begin{equation}
\Det_{\susn{2n}}(D)=\Deta(D_n)=\cP_{\psus{2n}}(D)\times_{\adet}\bbC.
\label{perdet.142}\end{equation}

\subsection{Isotropic determinant bundle}

One can make a different, but similar, construction using the isotropic
quantization of $D_{n}.$

\begin{definition}\label{detn.3}
For $\epsilon>0$, let $\De_{n}\in\Psi_{\wq{2n}{\epsilon}}^{1,1}(M/B;
2^{n-1}(E\oplus F))$ be the isotropic quantization of $D_{n}$ as in
Appendic~\ref{wq.0}, so giving an operator on $\cS(\bbR^{n}\times
X;2^{n-1}(E_{b}\oplus F_{b})).$
\end{definition}

As discussed earlier for families of standard elliptic operators, there are
two equivalent definitions of the determinant line bundle for $\De_{n}.$
Namely, Quillen's spectral definition or as an associated bundle to the
principal bundle of invertible perturbations. In the latter case, the principal
$G_{\wq{2n}{\epsilon}}^{-\infty}(M/B;2^{n-1}(E\oplus F))[[\varepsilon]]$-bundle
has fibre
\begin{multline}
\mathcal{P}_{\wq{2n}{\epsilon}}(\De_{en})_{b}=
\{\De_{2n,b}+Q_{b}; Q_{b}\in \Psi^{-\infty}_{\wq{2n}{\epsilon}}(
Z_{b}; 2^{n-1}(E_{b}\oplus F_{b})),\\
\exists\ (\De_{2n,b}+Q_{b})^{-1}\in 
\Psi_{\wq{2n}{\epsilon}}^{-1,-1}(Z_{b};2^{n-1}(E_{b}\oplus F_{b}))\}.
\label{detn.5}\end{multline}

Note that this fibre is non-empty as soon as the original family $D$ has
vanishing numerical index. Indeed, we know that $D_{n}$ then has vanishing
suspended index and hence has an invertible perturbation by a smoothing
operator (in the suspended sense). The isotropic product is smooth down to
$\epsilon =0,$ where it reduces to the suspended product (pointwise in the
parameters). Thus such a perturbation is invertible with respect to the
isotropic product for small $\epsilon >0.$ Since these products are all
isomorphic for $\epsilon >0,$ it follows that perturbations as required in
\eqref{detn.5} do exist.

\begin{proposition}\label{detdirac.32}
Let $D\in\Psi^{1}(M/B;E,F)$ be an elliptic family with vanishing numerical
index, then for each $n\in\bbN_{0}$ and $\epsilon>0,$ the determinant
line bundle $\Det(\De_{n})$ is naturally isomorphic to the determinant line
bundle $\Det(D)$.
\end{proposition}

\begin{proof}
The proof is by induction on $n\in \bbN_{0}$ starting with the trivial case
$n=0.$ We proceed to show that $\Det(\De_{n+1})\cong \Det(\De_{n}).$ In
Quillen's definition of the determinant line bundle, only the
eigenfunctions of the low eigenvalues are involved and the strategy is to
identify the eigensections of the low eigenvalues of $\De_{n}$ with those
of $\De_{n+1}.$ 

The isotropic quantization of the polynomial $\tau_n^2+t_n^2,$ is the
harmonic oscillator, $H_{n}^\epsilon,$ so
\begin{equation}
\De_{n+1,b}^{*}\De_{n+1,b}=
\begin{pmatrix}
                      \De_{n,b}^{*}\De_{n,b} +H_{n+1}^\epsilon-\epsilon & 0 \\
                      0 & \De_{n,b}\De_{n,b}^{*} +H_{n+1}^\epsilon+\epsilon   
                      \end{pmatrix}.
 \label{detdirac.44}\end{equation}
The eigenvalues of $H_{n+1}^\epsilon$ are positive, with the smallest being
simple. The eigensections of $\De^{*}_{n+1,b}\De_{n+1,b}$ and
$\De_{n+1,b}\De^{*}_{n+1,b}$ with smalle eigenvalues are of the form
\begin{equation}
\Phi^+(f_b)=
\begin{pmatrix}
\varphi_{n+1} \otimes f \\
0    
\end{pmatrix},\
\Phi^-(f_b)=\begin{pmatrix}
0 \\
\varphi_{n+1}\otimes \De_{n,b}f
\end{pmatrix}
\label{detdirac.46}\end{equation}
where $f$ is an eigenfunction of $\De^{*}_{n,b}\De_{n,b}$ with  
eigenvalue less than $2\epsilon.$  Note also that on such an eigenfunction, 
$\hat{D}_{n+1,b}$ acts as
\begin{equation}
\begin{pmatrix}
                     iC^{*}_{n+1} & \De_{n,b}^{*} \\
                      \De_{n,b} & iC_{n+1}   
\end{pmatrix}
\Phi^+(f_b)=\Phi^-(f_b)
\label{detdirac.47}\end{equation} 
since $C^{*}_{n+1}\varphi_{n+1}=0.$  For $0<\lambda<2\epsilon,$
consider the open set
\begin{equation}
\mathcal{U}_{\lambda}=\{b\in B; \lambda \, \,\mbox{is not an eigenvalue of} 
\, \,D_{b}^{*}D_{b}\}. 
\label{detdirac.49}\end{equation}
Let $\cH^{+,k}_{[0,\lambda)}$ denote the vector bundle over
$\mathcal{U}_{\lambda}$ spanned by the eigenfunctions of
$\De_{k,b}^{*}\De_{k,b}$ with eigenvalues less than $\lambda$.  Let
$\cH^{-,k}_{[0,\lambda)}$ denote the vector bundle over
$\mathcal{U}_{\lambda}$ spanned by the eigenfunctions of $\De_{k,b}
\De_{k,b}^{*}$ with eigenvalues less than $\lambda.$
Then there are natural identifications
\begin{equation}
 F_{\mathcal{U},\lambda}^{\pm,n}:\cH^{\pm,n}_{[0,\lambda)}\ni f_b
\longmapsto \Phi^\pm(f_b)\in\cH^{\pm,n+1}_{[0,\lambda)}.
\label{detdirac.48}\end{equation}
Thus, directly from Quillen's definition of the determinant bundle $\De_{n+1,b}$ 
and $\De_{n,b}$ have isomorphic determinant line bundles.
\end{proof}

\subsection{Adiabatic limit of $\Det(\De_{n})$}

\begin{proposition}\label{perdet.143} If
  $P\in\Psi^{m}_{\psusn{2n}}(M/B;E,F)$ is a family of fully elliptic
  operators with vanishing numerical index then the bundle over
  $B\times[0,1]$ with fibre 
\begin{multline}
\cP_{b,\epsilon }=\big\{Q\in\Psi^{-\infty}_{\susn{2n}}(Z_b,E_b);\\
\exists\ (P+Q)^{-1}\in\Psi^{-m}_{\psusn{2n}}(Z_b;F_b,E_b)
\text{ for the }\epsilon\text{ isotropic product}.\big\}
\label{perdet.144}\end{multline}
is a principal $\cG$-bundle for the bundle of groups with fibre
\begin{multline}
\cG^{-\infty}_{\susn{2n},\epsilon}(Z_b,E_b)=
\big\{\Id+A,\ A\in\Psi^{-\infty}_{\susn{2n}}(Z_b,E_b);\\
\exists\ (\Id+A)^{-1}=\Id+B,\ B
\in\Psi^{-\infty}_{\susn{2n}}(Z_b,E_b)
\text{ for the }\epsilon\text{ isotropic product}\big\}
\label{perdet.145}\end{multline}
and the associated determinant bundle defined over $\epsilon >0$ extends
smoothly down to $\epsilon =0$ and at $\epsilon =0$ is induced by the
adiabatic determinant.
\end{proposition}

\begin{proof} This is just the smoothness of the `rescaled' determinant
  (\ie with the singular terms removed) down to $\epsilon =0.$
\end{proof}

We will now complete the proof of Theorem~\ref{int.6} in the Introduction
which we slightly restate as

\begin{theorem}[Periodicity of the determinant line bundle] Let $D\in
\Psi^{1}(M/B;E,F)$ be an elliptic family with vanishing numerical index,
then for $n\in \bbN$ and $\epsilon>0,$
\begin{equation*}
   \Deta(D_{n})\cong \Det(\De_{n})\cong \Det(D).
\end{equation*}
\label{detnp.7}\end{theorem}

\begin{proof} The existence of the second isomorphism follows from
Proposition~\ref{detdirac.32}. The first follows from
Proposition~\ref{perdet.143}.
\end{proof}  

\section{Eta invariant}
\label{eta.0}

In \cite{MR96h:58169} a form of the eta invariant was discussed for
elliptic and invertible once-suspended families of pseudodifferential
operators. Applied to the spectral family (on the imaginary axis) of a
self-adjoint invertible Dirac operator this new definition was shown to
reduce to the original definition, of Atiyah, Patodi and Singer in
\cite{MR53:1655b} of the eta invariant of a single operator. Here, the
definition in \cite{MR96h:58169} is shown to extend to (fully) elliptic,
invertible, product-suspended families. In the \S\ref{Adeta} it is further
extended to such product-suspended families in any odd number of variables. The
extension to single-parameter product-suspended operators allows us to
apply the definition to $A+i\tau,$ $\tau\in\bbR,$ for $A\in\Psi^1(X;E)$ an
invertible elliptic selfadjoint pseudodifferential operator and check that
this reduces to the spectral definition, now as given by Wodzicki
(\cite{Wodzicki3}). Again the extended (and below also the `adiabatic') eta
invariant gives a log-multiplicative function for invertible families
\begin{equation}
\eta(AB)=\eta(A)+\eta(B)
\label{perdet.40}\end{equation}
and this allows us to show quite directly that the associated $\tau$ invariant
is a determinant in the sense discussed above.

\subsection{Product-suspended eta}

If $B\in\Psi_{\psus}^{m,m'}(X;E)$ is a product-suspended family it satisfies 
\begin{equation}
\frac{\pa^N}{\pa\tau^N}B(\tau)\in \Psi_{\psus}^{m-N,m'-N}(X;E)\ \forall\
N\in\bbN_0.
\label{perdet.41}\end{equation}
This implies that for $N$ large, say $N>\dim X+m,$ the differentiated
family takes values in operators of trace class on $L^2.$

\begin{proposition}\label{perdet.42} For any $m,$ $m'\in\bbZ,$ if $N\in\bbN$ is
chosen sufficiently large then,
\begin{equation}
B\in\Psi_{\psus}^{m,m'}(X;E)\Longrightarrow
\Tr_E\left(\frac{\pa^N}{\pa\tau^N}B(\tau)\right)\in\CI(\bbR^{p})
\label{perdet.43}\end{equation}
has a complete asymptotic expansion (possibly with logarithms) as
$\tau\to\pm\infty$ and the coefficient of $T^0$ in the expansion as
$T\to\infty$
\begin{equation}
\begin{gathered}
\bTr(B)=\LIM_{T\to\infty}F_{B,N}(T),\\
F_{B,N}(T)=\int_{-T}^T\int_0^{t_1}\dots\int_0^{t_N}
\Tr_E\left(\frac{\pa^N}{\pa s^N}B(s)\right)
ds\, dt_N\dots dt_1
\end{gathered}
\label{perdet.44}\end{equation}
is independent of the choice of $N$ and defines a trace functional 
\begin{equation}
\bTr:\Psi_{\psus}^{\bbZ,\bbZ}(X;E)\longrightarrow \bbC,\
\bTr([A,B])=0\ \forall\ A,B\in\Psi_{\psus}^{\bbZ,\bbZ}(X;E)
\label{perdet.45}\end{equation}
which reduces to 
\begin{equation}
\bTr(B)=\int_{\bbR}\Tr_E\left(B(\tau)\right) d\tau\ \forall\
A\in\Psi_{\psus}^{-\infty,-\infty}(X;E).
\label{perdet.61}\end{equation}
\end{proposition}

\begin{proof} As already noted, $\pa_s^NB(s)$ is a continuous family
of trace class operators as soon as $N>\dim X+m.$ Then
\eqref{perdet.43} is a continuous function and further differentiation
again gives a continuous family of trace class operators so the trace is
smooth.

To see that this function has a complete asymptotic expansion we
appeal to the discussion of the structure of the kernels of such
product-suspended families in Appendix~\ref{iso.0}. It suffices to consider the
trace of a general element $B\in\Psi^{-n-1,0}_{\psus}(X;E).$ Since
the kernels form a module over $\CI(X^2)$ we can localize in the base
variable (not directly in the 
suspended variable since that has global properties). Localizing near a
point away from the diagonal gives a classical symbol in the suspending
variable with values in the smoothing operators. Since the trace is the
integral over the diagonal this makes no contribution to
\eqref{perdet.43}. Thus it suffices to suppose that $B$ is supported in the
product of a coordinate neighbourhood with itself over which the bundle $E$
is trivial. Locally (see \eqref{fipomb2.7}) the kernel is given by Weyl
quantization of a product-type symbol so the trace becomes the integral of
the sum of the diagonal terms and hence we need only consider
\begin{equation}
\frac{1}{(2\pi)^n}\int a(x,\xi,\tau)dxd\xi
\label{perdet.46}\end{equation}
where $a$ is compactly supported in the base variables $x.$ Now by
definition, $a$ is a smooth function, with compact support, on the product
$\bbR^n\times[\com{\bbR\times\bbR^n};\pa(\com{\bbR}\times\{0\})].$ Thus
we can further localize the support of $a$ on this blown up space. There
are three essentially different regions, corresponding to the part of the
boundary which arises from the radial compactfication, the part arising from the
blow up and the corner.

The first of these regions corresponds to a true suspended family, as
considered in \cite{MR96h:58169}. In this region the variable $|\xi|$
dominates, and $|\tau|\le C|\xi|$ on the support so the integral takes the
form
\begin{multline}
\frac{1}{(2\pi)^n}\int_0^1\int \phi(r\tau)
r^{n+1}f(x,\omega ,r,r\tau)r^{-n-1}dxd\omega dr\\
=\int_{0}^{1}\int s\phi (R)f(x,\omega ,Rs,R)dxd\omega dR,\ s=1/\tau.
\label{perdet.47}\end{multline}
Here, $\phi$ has compact support and $f$ (with the factor of $r^{n+1}$
representing the order $-n-1)$ is smooth. The result is smooth in
$s=1/\tau,$ which corresponds to a complete asymptotic expansion with only
non-negative terms.

The second region corresponds to boundedness of the variable $\xi$ with the
function being a classical symbol (by assumption of order at most $0)$ in
$\tau$ so integration simply gives a symbol 
\begin{equation}
\frac{1}{(2\pi)^n}\int a(x,\xi,\tau)dxd\xi.
\label{perdet.48}\end{equation}

The third region is the most problematic. Here the two boundary
faces of the compactification are defined by $r=1/|\xi|$ and $|\xi|/\tau$ and
with polar variables $\omega =\xi/|\xi|.$ Thus the integral takes the form 
\begin{equation}
\int r^{n+1}f(x,\omega,r,s/r) r^{-n-1}dxd\omega dr\in\CI([0,1)_s)+(\log
  s)\CI([0,1)_s),
\label{perdet.49}\end{equation}
where $f$ is smooth and with compact support near $0$ in the last two
variables. This a simple example of the general theorem on pushforward
under b-fibrations in \cite{MR93i:58148}, or the `singular asymptotics lemma' of
Br\"uning and Seeley (see also \cite{Grieser-Gruber1}) and is in fact a
type of integral long studied as an orbit integral. In any case the
indicated regularity follows and this proves the existence of a complete
asymptotic expansion, possibly with single logarithmic terms.

It follows that the integral in \eqref{perdet.44} also has a complete
asymptotic expansion as $T\to\infty;$ where in principle there can be
factors of $(\log T)^2$ after such integration. Thus the
coefficient of $T^0$ does exist, and defines $\bTr(B).$ Now if $N$ is
increased by one in the definition, the additional integral gives the same
formula \eqref{perdet.44} except that a constant of integration may be
added by the first integral. After $N$ additional integrals, this adds a
polynomial, so the result is changed by the integral over
$[-T,T]$ of a polynomial. This is an odd polynomial, so has no constant
term in its expansion at infinity. Thus the definition of $\Tr(B)$ is in
fact independent of the choice of $N.$

The trace identity follows directly from \eqref{perdet.44}, since if
$B=[B_1,B_2]$, then any derivative is a sum of commutators between operators with
order summing to less than $-n$ and the trace of such a term vanishes. Thus
applied to a commutator \eqref{perdet.44} itself vanishes.
\end{proof}

Using this trace functional on product-suspended operators we extend the
domain of the eta invariant.

\begin{proposition}\label{perdet.53} The eta invariant for any
fully elliptic, invertible element $A\in\Psi^{m,m}_{\psus}(X;E)$
defined using the regularized trace
\begin{equation}
\eta(A)=\frac1{\pi i}\bTr(A^{-1}\dot A),\ \dot A=\frac{\pa A}{\pa\tau}
\label{perdet.50}\end{equation}
is a log-multiplicative functional, in the sense of \eqref{perdet.40}.
\end{proposition}

\begin{proof} Certainly \eqref{perdet.50} defines a continuous functional
on elliptic and invertible product\--suspended families. The log-multiplicativity
formula \eqref{perdet.40} follows directly since if $B$ is another
invertible product-suspended family then
\begin{equation}
(AB)^{-1}\frac{\pa(AB)}{\pa\tau}=B^{-1}A^{-1}\dot A B+B^{-1}\dot B
\label{perdet.51}\end{equation}
and the trace identity shows that $\bTr(B^{-1}A^{-1}\dot A B)=\eta(A).$
\end{proof}

\subsection{$\eta(A+i\tau)=\eta (A)$}

To relate this functional on product-suspended invertible operators to the
more familiar eta invariant for self-adjoint elliptic pseudodifferential
operators we rewrite the definition in a form closer to traditional zeta
regularization, starting with the regularized trace.

Consider the meromorphic family $t^{-z}_+$ of tempered
  distributions with support in $[0,\infty).$ This family has poles only at
  the positive integers, with residues being derivatives of the delta
  function at the origin. For $\Re z$ sufficiently positive and
  non-integral, $t^{-z}_+$ can be paired with the function $F_{B,N}(t)$ in
  \eqref{perdet.44}, since this is smooth and of finite growth at infinity.
  This pairing gives a meromorphic function in $\Re z>C,$ with poles only at the
  natural numbers since the poles of $t^{-z}_+$ are associated with the
  behaviour at $0,$ where $F_{B,N}$ is smooth. In fact this pairing 
\begin{equation}
g(z)=\langle T^{-z-1}_+,F_{B,N}(T)\rangle
\label{perdet.55}\end{equation}
extends to be meromorphic in the whole complex plane. Indeed, dividing the
  pairing into two using a cut-off $\psi\in\CIc([0,\infty))$ which is
  identically equal to $1$ near $0,$  
\begin{equation}
g(z)=\langle T^{-z-1}_+,\psi (T)F_{B,N}(T)\rangle+
\langle T^{-z-1}_+,(1-\psi (T)))F_{B,N}(T)\rangle 
\label{perdet.56}\end{equation}
the first term is meromorphic with poles only at $z\in\bbN$ and the poles
  of the second term arise from the terms in the asymptotic expansion of
  $F_{B,N}(T).$ Notice that there is no pole at $z=0$ for the first term
  since the residue of $T^{-z-1}_+$ at $z=0$ is a multiple of the delta
  function and $F_{B,N}(0)=0.$ The pole at $z=0$ for the second term arises
  exactly from the coefficient of $T^0$ in the asymptotic expansion so we
  see that 
\begin{equation}
\bTr(B)=\res_{z=0}g(z).
\label{perdet.64}\end{equation}
Any terms $a_k(\log T)^k$ for $k\in\bbN,$ in the expansion do not
  contribute to the residue since they integrate to regular functions at
  $z=0$ plus multiples of $z^{-k}.$

\begin{proposition}\label{perdet.54} For $B\in\Psi^{m,m'}_{\psus}(X;E)$
  and any $N>m-\dim X-1,$ the regularized trace is the residue
  at $z=0$ of the meromorphic continuation from $\Re z>>0,$ $z\notin\bbZ,$ of
\begin{equation}
\frac{(-1)^{N+1}}{(N-z)\dots(1-z)(-z)}
\left\langle \left(\frac{(t+i0)^{N-z}}{1+e^{-\pi i z}}+
\frac{(t-i0)^{N-z}}{1+e^{\pi i z}}\right),\Tr_E(\pa_t^N B(t)\right\rangle.
\label{perdet.57}\end{equation}
\end{proposition}

\begin{proof} Consider the identity 
\begin{equation}
t^{-z-1}_+=\frac1{(N-z)\dots(1-z)(-z)}\frac{d^{N+1}}{dt^{N+1}}t^{-z+N}_+.
\label{perdet.58}\end{equation}
After inserting this into \eqref{perdet.56}, integration by parts is
justified (since \eqref{perdet.58} holds in the sense of distributions on
the whole real line, supported in $[0,\infty)),$ and shows that 
\begin{multline}
g(z)=\\
\frac1{(N-z)\dots(1-z)(-z)}\langle t^{-z+N}_+,
(-1)^{N+1}\Tr_E(\frac{\pa^N}{\pa t^N}B)(t)-\Tr_E(\frac{\pa^N}{\pa t^N}B)(-t)
\rangle
\label{perdet.59}\end{multline}
where the pairings are defined, and holomorphic, for $\Re z$ large and $z$
non-integral. Using the identity 
\[
     t_+^z=\frac{(t+i0)^z}{1-e^{2\pi i z}} +\frac{(t-i0)^z}{1-e^{-2\pi i z}},
\]
\eqref{perdet.59} becomes 
\begin{equation}
\begin{gathered}
g(z)=\frac{(-1)^{N+1}}{(N-z)\ldots(1-z)(-z)}
\langle D(t,z),\Tr_E(\frac{\pa^N}{\pa t^N}B)(t)\rangle,\\
{\begin{aligned}
D(t,z)=\frac{(t+i0)^{N-z}}{1-e^{2\pi i(N-z)}}+&
\frac{(t-i0)^{N-z}}{1-e^{-2\pi i(N-z)}}\\
&+ (-1)^{N}\frac{(-t+i0)^{N-z}}{1-e^{2\pi i(N-z)}}+
(-1)^{N}\frac{(-t-i0)^{N-z}}{1-e^{-2\pi i(N-z)}}.
\end{aligned}}
\end{gathered}
\label{perdet.72}\end{equation}
Now, $(-t-i0)^{-z}=e^{\pi iz}(t+i0)^{-z}$ so
\begin{multline}
D(t,z)=\left(\frac{1}{1-e^{-2\pi i z}}
+\frac{e^{\pi i z}}{1-e^{2\pi i z}}\right)
(t+i0)^{N-z} \\
+\left(\frac{1}{1-e^{2\pi i z}}+\frac{e^{-\pi i z}}{1-e^{-2\pi i z}}\right)
(t-i0)^{N-z}
\label{perdet.74}\end{multline}
which reduces \eqref{perdet.72} to \eqref{perdet.57}.
\end{proof}

This allows us to prove a result of which Theorem~\ref{int.7} in the
introduction is an immediate corollary.

\begin{theorem}\label{perdet.52} If $A\in\Psi^1(X;E)$ is a self-adjoint
  elliptic and invertible pseudodifferential operator then $\eta(A+i\tau),$
  defined through \eqref{perdet.50} reduces to the (regularized) value at
  $z=0$ of the analytic continuation from $\Re z>>0$ of
\begin{equation}
\sum\limits_{j}\sgn(a_j)|a_j|^{-z},
\label{perdet.66}\end{equation}
where the $a_j$ are the eigenvalues of $A,$ in order of increasing $|a_j|$
  repeated with multiplicities.
\end{theorem}

\begin{proof} With $A(\tau)=A+i\tau$ the eta invariant defined by
\eqref{perdet.50} reduces to 
\begin{equation}
\eta(A+i\tau)=\frac1{\pi}\bTr\left((A+i\tau)^{-1}\right)=\frac1\pi\res_{z=0}h(z)
\label{perdet.67}\end{equation}
where $h(z)$ is the function \eqref{perdet.57} with $B(t)=(A+it)^{-1}.$
Computing the $N$th derivative
\begin{equation}
\frac{\pa^N}{\pa\tau^N}\left((A+i\tau)^{-1}\right)
=i(-1)^{N+1}N!(\tau-iA)^{-N-1}.
\label{perdet.68}\end{equation}
The trace is therefore given, for any $N>n,$ by 
\begin{equation}
\tr_E\left(\frac{\pa^N}{\pa\tau^N}(A+i\tau)^{-1}\right)=i(-1)^{N+1}N!
\sum\limits_{j}(\tau-ia_j)^{-N-1}.
\label{perdet.69}\end{equation}
This converges uniformly with its derivatives so can be inserted in the
pairing \eqref{perdet.57} and the order exchanged. Thus  
\begin{multline}
h(z)=\frac{a_N(z)i(-1)^{N+1}N!}{1+e^{-\pi i z}}\sum\limits_{j}
\lim_{\epsilon \downarrow0}\int_{\bbR+i\epsilon}
\tau^{N-z}(\tau-ia_j)^{-N-1}d\tau
\\
+\frac{a_N(z)i(-1)^{N+1}N!}{1+e^{\pi i z}}\sum\limits_{j}\lim_{\epsilon
  \downarrow0}\int_{\bbR-i\epsilon}\tau^{N-z}(\tau-ia_j)^{-N-1}d\tau, 
\label{perdet.70}\end{multline}
where 
\begin{equation*}
      a_{N}(z)= \frac{(-1)^{N+1}}{(N-z)\ldots (1-z)(-z)}.
\label{perdet.163}\end{equation*}
Each of these integrals are contour integrals, actually independent of
$\epsilon >0$ for $\epsilon$ smaller than the minimal $|a_j|.$ By residue
computation, in the first sum by moving the contour to infinity in the
upper half plane and in the second by moving the contour into the lower
half plane
\begin{equation}
\int_{\bbR\pm i\epsilon}\tau^{N-z}(\tau-ia_j)^{-N-1}d\tau=
\begin{cases}\pm2\pi i \frac{(N-z)\cdots (1-z)}{N!}e^{\mp\pi i z/2}|a_j|^{-z}&
\pm a_j>0
\\
0&\pm a_j<0
\end{cases}
\label{perdet.71}\end{equation}
Inserting this into \eqref{perdet.70} shows that $\eta(A+i\tau)$ is the
residue at $z=0$ of
\begin{equation}
\frac1{z\cos(\pi z/2)}\sum\limits_{j}\sgn(a_j)|a_j|^{-z}.
\label{perdet.65}\end{equation}
By definition, the
usual eta invariant, $\eta(A),$ is the value at $z=0$ of the continuation of
the series in \eqref{perdet.65}. This series is the analytic continuation
of the trace of an entire family of classical elliptic operators of complex
order $-z$ (namely $A^{-z}(\Pi_+-\Pi_-)$ where $\Pi_{\pm}$ are the
projections onto the span of positive and negative eigenvalues) which can
have only a simple pole at $z=0.$ In fact, here, it is known that there is no
singularity, \ie the residue vanishes. Even without invoking this we
conclude the desired equality, since the explicit meromorphic factor in
\eqref{perdet.65} is odd in $z,$ so a pole in the continuation of the
series would not affect the residue.
\end{proof}

\section{Universal $\eta$ and $\tau,$ invariants}\label{etat.0}

That the differential of the eta invariant of a family of self-adjoint
Dirac operators is a multiple of the first (odd) Chern class of the index,
in odd cohomology of the base, of the family is well-known. In the case of
the suspended eta invariant discussed in \cite{MR96h:58169} and above, we
show that the $\eta$ invariant is, in appropriate circumstances, the
logarithm of a determinant, which is to say a multiplicative function
giving the first odd Chern class. Initially we show this in the context of
classifiying spaces for K-theory, then in the geometric context of
$(2n+1)$-fold suspended odd elliptic families.

Consider again the algebra of once-suspended isotropic pseudodifferential
operators of order $0$ on $\bbR^n,$ with values in smoothing operators on a
compact manifold $X.$ This can be identified with the smooth functions on
$\com{\bbR^{2n+1}}\times X^2$ and the subspace 
\begin{equation}
\cI_+=\left\{A\in\CI(\com{\bbR^{2n+1}}\times X^2);
A\cong0\Min\{t\le0\}\cap\bbS^{2n}\times X^2\right\},
\label{perdet.75}\end{equation}
is a subalgrebra. Here, $t$ is the suspending parameter and equality is in 
the sense of Taylor series at infinity on the compactified Euclidean
space. Thus the subalgebra is just the sum of the 
smoothing ideal (identified with the functions vanishing to infinite
order everywhere at the boundary) and the subalgebra of functions vanishing
in $t<0.$ In fact $\cI_+$ is also an ideal. We consider the corresponding group 
\begin{equation}
\cG_+=\{B=\Id+A,\ A\in\cI_+,\ B^{-1}=\Id+B',\ B'\in\cI_+\}.
\label{perdet.76}\end{equation}

Now we may use the suspending variable $t$ to identify the upper
half-sphere $\{t>0\}\cap\bbS^{2n}$ of the boundary of $\com{\bbR^{2n+1}}$
with $\bbR^{2n},$  
\begin{equation}
\{t>0\}\cap\bbS^{2n}\ni[(t,x,\xi)]\longmapsto (X,\Xi)=(x/t,\xi/t)\in\bbR^{2n}.
\label{perdet.77}\end{equation}
The inverse image under pull-back of $\cS(\bbR^{2n})$ is then naturally
identified with $\{a\in\CI(\bbS^{2n});a=0\Min t<0\}$ where $\bbS^{2n}$ is
the boundary of the radial compactification of $\bbR^{2n+1}.$
This allows the space of formal power series $\cS(\bbR^{2n})[[t]]$
to be identified with the formal power series at the boundary of the
subspace of $\CI(\com{\bbR^{2n+1}})$ consisting of the functions vanishing
in $t<0.$

The same identifications carry over to the case of functions valued in the
smoothing operators and so gives a short exact sequence of algebras
\begin{equation}
\xymatrix{\Psi^{-\infty}_{\sus}(\bbR^n\times X)\ar[r]&
\cI_+
\ar[r]&
\cS(\bbR^{2n}\times X^2)[[t]]}
\label{perdet.78}\end{equation}

\begin{lemma}\label{perdet.79} In \eqref{perdet.78}, the product induced on
  the quotient is the standard $\star$ product (valued in smoothing
  operators on $X)$ on $\bbR^{2n}$ (\ie the `Moyal product').
\end{lemma}

\begin{proof} 
Let $A, B\in \cI_+$ be the symbols of two operators $\hat{A}, \hat{B}\in
\Psi^{0}_{\psus(1)}(\bbR^n\times X)$. Then the asymptotic expansion at
infinity of the symbol of $\hat{A}\hat{B}$ is given by the standard $\star$
product
\begin{equation}
  \left. \sigma(\hat{A}\hat{B})\sim \sum_{k=0}^{\infty}  \frac{1}{k!(2i)^{k}}
   (D_{x}D_{\eta}- D_{y}D_{\xi})^{k}A(t,x,\xi)B(t,y,\eta) 
\right|_{x=y, \eta=\xi}.
\label{star.1}\end{equation} 
Under the map \eqref{perdet.77}, 
the asymptotic expansion \eqref{star.1} becomes
an asymptotic expansion at $\{t>0\}\cap\bbS^{2n}\subset\com{\bbR^{2n+1}}$
\begin{equation}
\left. \sigma(\hat{A}\hat{B})\sim \sum_{k=0}^{\infty}  \frac{1}{k!(2i)^{k}}
  \frac{1}{t^{2k}}
   (D_{X}D_{\Lambda}- D_{Y}D_{\Xi})^{k}A(t,tX,t\Xi)B(t,tY,t\Lambda) 
\right|_{X=Y, \Lambda=\Xi}.
\label{star.2}\end{equation}
Thus, if 
\begin{equation}
A(t,tX,t\Xi) \sim \sum_{k}^{\infty} \frac{1}{t^{k}} a_{k}(X,\Xi),\
B(t,tX,t\Xi) \sim \sum_{k}^{\infty} \frac{1}{t^{k}} b_{k}(X,\Xi)
\label{star.3}\end{equation}
are the asymptotic expansions of $A$ and $B$ at
$\{t>0\}\cap\bbS^{2n}\subset\com{\bbR^{2n+1}},$ then
\begin{equation}
\sigma(\hat{A}\hat{B}) \sim \sum_{k,l,m\ge 0} \frac{1}{k!(2i)^{k}}
\frac{1}{t^{2k+l+m}} \left.
(D_{X}D_{\Lambda}- D_{Y}D_{\Xi})^{k}a_{l}(X,\Xi)b_{m}(Y,\Lambda) 
\right|_{X=Y, \Lambda=\Xi}
\label{star.4}\end{equation}
is the asymptotic expansion of $\sigma(\hat{A}\hat{B})$ at 
$\{t>0\}\cap\bbS^{2n}\subset\com{\bbR^{2n+1}}$.  But the right hand side
is precisely the standard $\star$ product on 
$\cS(\bbR^{2n}\times X^2)[[\varepsilon]]$ with $\varepsilon=\frac{1}{t^{2}}$.
\end{proof}

Corresponding to this exact sequence of algebras is the exact sequence
of groups consisting of the invertible perturbations of the identity 
\begin{equation}
G^{-\infty}_{\sus}(\bbR^n\times X)\longrightarrow\cG_+\longrightarrow
G^{-\infty}(\bbR^n\times X)[[t]].
\label{perdet.80}\end{equation}

\begin{theorem}\label{perdet.81} In this `delooping' sequence, the
first group is classifying for even K-theory, the central group is (weakly)
contractible and the quotient is (therefore) a classifying group for odd
K-theory; the eta invariant, defined as in \eqref{perdet.50},
\begin{equation}
\eta:\cG_+\longrightarrow \bbC
\label{perdet.82}\end{equation}
restricts to twice the index on the normal subgroup and $e^{i\pi\eta}=\deta$ is
the adiabatic determinant on $G^{-\infty}(\bbR^n\times X)[[t]].$
\end{theorem}

\begin{proof} As a first step in the proof we consider the behaviour of the
  regularized trace.

\begin{lemma}\label{perdet.83} The regularized trace $\bTr$ on the central
  algebra in \eqref{perdet.78} restricts to the integrated trace on the
  smoothing subalgebra and 
\begin{equation}
\bTr(\frac{\pa b}{\pa t})=\int_{\bbR^{2n}} b_{2n}dXd\Xi
\label{perdet.84}\end{equation}
for any $b\in \cI_+,$ where $b_k$ is the term of order $k$ in the formal
  power series of the image in \eqref{perdet.78}.
\end{lemma}

\begin{proof} When the parameter $t$ is fixed, an element $b\in\cI_+$ is
  actually a smoothing operator, since the asymptotic behavior on the
  surface where $t$ is constant is determined by the equatorial sphere
  $t=0$ at infinity. Thus the definition, from \eqref{perdet.44}, of
  $\bTr(b)$ for any element $b\in\cI_+$ may be modified by dropping all $N$
  integrals, \ie we may take $N=0.$ Indeed, taking $N>0$ and then
  integrating results in the case $N=0,$ plus a polynomial which, as noted
  earlier, does not affect the result. Carrying out the last integral
  by the fundamental theorem of calculus,
\begin{equation}
\bTr(\dot b)=\LIM_{T\to\infty}
(\int_{\bbR^{2n}}b(T,x,\xi) dxd\xi-\int_{\bbR^{2n}} b(-T,x,\xi) dxd\xi)
\label{perdet.85}\end{equation}
where $\LIM$ stands for the constant term in the asymptotic expansion.
The second term in \eqref{perdet.85} corresponds to $t<0$ where $b$ is
rapidly decreasing so does not contribute to the asymptotic expansion.
Now, making the scaling change of variable in \eqref{perdet.77}, transforms
\eqref{perdet.85} to
\begin{equation}
\bTr(\dot b)=\LIM_{T\to\infty}T^{-2n}\int_{\bbR^{2n}}\tilde b(T,X,\Xi)dXd\Xi
\label{perdet.86}\end{equation}
where $\tilde b$ is the transformed function. Thus \eqref{perdet.86} picks
out the term of homogeneity $2n$ (in $T)$ in the formal expansion of
$\tilde b.$ This gives exactly \eqref{perdet.84}.
\end{proof}

Now, by definition, the eta invariant is $\frac1{\pi i}\bTr(a^{-1}\dot a).$
It follows directly that restricted to the smoothing subgroup this lies in
$2\bbZ.$ Thus $D=e^{i\pi\eta}$ does indeed descend to the quotient group in
\eqref{perdet.80}. This group is connected, so to check that it reduces to
the `adiabatic' determinant defined earlier we only need check the
variation formula, both being $1$ on the identity. Along a curve $a(s),$ 
\begin{equation}
\frac{d}{ds}\eta(a(s))=\frac1{\pi i}\bTr(a^{-1}\frac{d\dot
  a}{ds}-a^{-1}\frac{da}{ds} a^{-1}\dot a)
=\bTr\left(\frac{d}{dt}(a^{-1}\frac{da}{ds})\right).
\label{perdet.87}\end{equation}
Thus the identity \eqref{perdet.84} shows that 
\begin{equation}
\frac{d}{ds}\eta(a(s))=\Tr\left[(\tilde a(s)\frac{d\tilde a}{ds})_{2n}\right]
\label{perdet.88}\end{equation}
where $\tilde a$ is the image of $a$ in the third group in
\eqref{perdet.80}. The identity term in $a$ does not affect the argument
since it is annihilated by $d/ds.$

Since the right hand side of \eqref{perdet.88} is the variation formula for
the logarithm of the adiabatic determinant this proves the theorem.
\end{proof}

\section{Geometric $\eta$ and $\tau$ invariants}\label{Geometa}

Returning to the `geometric setting' of a fibration \eqref{int.4} with
compact fibres, consider a totally elliptic family
$A\in\Psi^{m,m'}_{\psus}(M/B;E,F).$ Although we allow for operators between
different bundles here, \eqref{perdet.50} is still meaningful as a
definition of the eta invariant if $A$ is invertible. Consider the
principal bundle, of the type discussed above,
\begin{equation}
\xymatrix{G^{-\infty}_{\sus}(M/B;E)\ar@{-}[r]&\mathcal{A}\ar[d]^{\nu}\\&B}
\label{perdet.166}\end{equation}
with fibre 
\begin{equation}
\mathcal{A}_b=\left\{A+Q;Q\in\Psi^{-\infty}_{\sus}(Z_b;E_b,F_b),\
(A+Q)^{-1}\in\Psi^{-m,-m'}_{\psus}(Z_b;F_b,E_b)\right\}.
\label{perdet.167}\end{equation}

\begin{proposition}\label{perdet.168} The eta invariant, defined by
  \eqref{perdet.50}, is a smooth function on $\cA$ such that for the fibre
  action of the structure group at each point
\begin{equation}
\eta(A(\Id+L))=\eta(A)+2\ind(\Id+L)
\label{perdet.169}\end{equation}
so projects to
\begin{equation}
\tau=e^{i\pi\eta }:B\longrightarrow \bbC^*
\label{perdet.170}\end{equation}
which represents the first odd Chern class of the index bundle of the
  family $A.$
\end{proposition}
\noindent In particular this result applies to an elliptic, self-adjoint,
  family of pseudodifferential operators of order $1$ by considering the
  spectral family.

\begin{proof} That $\eta:\cA\longrightarrow \bbC$ is well-defined follows
  from the discussion above as does the multiplicativity
  \eqref{perdet.169}. Thus, $\tau$ is well-defined as a function on $B$ and
  it only remains to check the topological interpretation.

Note that the fibre of $\cA$ is non-empty at each point of the
base. In fact it is always possible to find a global smoothing perturbation
to make the family invertible, although only when the families index
vanishes is this possible with a smoothing perturbation of compact support in the
parameter space. Thus, in complete generality, it is possible to choose a
smooth map 
\begin{multline}
Q_+:\bbR\longrightarrow \Psi^{-\infty}(M/B;E,F)\text{ such that}\\
Q_+(t)=0\Mfor t<<0,\ Q_+(t)=Q_+(T)\Mfor t\ge T>>0,\\
(A(t)+Q(t))^{-1}\in\Psi^{-m,-m'}(M/B;E,F)\ \forall\ t\in\bbR.
\label{perdet.171}\end{multline}
This follows directly from the fact that the index bundle, over $\bbR\times
B,$ is trivial for $t<<0$ and so is trivial over $\bbR\times B$ but
defines a generally non-trivial index class in $K^1(B).$ In fact the index
class of the family $A$ is represented by the map 
\begin{equation}
(A(T)+Q_+(T))^{-1}A(T)\in G^{-\infty}(M/B;E).
\label{perdet.172}\end{equation}
Namely, if this family is deformable to the identity in this bundle of
groups then there is a perturbation of compact support in $t$ making the
original family invertible.

The existence of $Q_+$ may be directly related to a larger principal bundle
with bundle of structure groups $\cG^{-\infty}_{\sus,+}(M/B;E)$ with fibre
\begin{multline}
\big\{\Id+Q_;\Id+Q_b\in\CI(\bbR;G^{-\infty}(Z_b;E_b),
\rho(t)Q(t)\in\cS(\bbR;\Psi^{-\infty}(Z_b;E_b),\\ \exists\
Q_0\in\Psi^{-\infty}(Z_b;E_b),\ (1-\rho
(t))(Q(t)-Q_0)\in\cS(\bbR;\Psi^{-\infty}(Z_b;E_b))\big\}.
\label{perdet.173}\end{multline}
Here $\rho (t)\in\CI(\bbR)$ is equal to $1$ in $t<-1$ and vanishes in
$t>1.$ Thus the short exact sequence of groups 
\begin{equation}
G^{-\infty}_{\sus}(M/B;E)\longrightarrow
G^{-\infty}_{+,\sus}(M/B;E)\overset{\pi_\infty}\longrightarrow
G^{-\infty}(M/B;E)
\label{perdet.174}\end{equation}
is the `delooping sequence' for $G^{-\infty}(M/B;E).$ In particular the
central group is weakly contractible and we may consider the enlarged
principal bundle 
\begin{equation}
\xymatrix{G^{-\infty}_{+,\sus}(M/B;E)\ar@{-}[r]&\mathcal{A}_+\ar[d]\\&B}
\label{perdet.175}\end{equation}
defined by replacing $G^{-\infty}_{\sus}$ above by $G^{-\infty}_{+,\sus}.$

The existence of $Q_+$ shows that this bundle is trivial, \ie has a global
section  
\begin{equation*}
q:B\longrightarrow \cA_+
\label{perdet.176}\end{equation*}
which induces a `classifying bundle map'
\begin{equation*}
\tq:\cA\longrightarrow \cG^{-\infty}_{+,\sus}(M/B;E),
\tq(A_b+Q_b)=(A_b+Q_{+,b})^{-1}(A_b+Q_b)\in \cG^{-\infty}_{+,\sus}(Z_b;E_b).
\label{perdet.178}\end{equation*}

Now, the definition and basic properties of the eta invariant given by
\eqref{perdet.50} are quite insensitive to the enlargment of $\cA$ to
$\cA_+$ and so still define a smooth function $\eta^+:\cA_+\longrightarrow
\bbC.$ The same is true for the group $\cG^{-\infty}_{+,\sus}(M/B;E),$
defining the corresponding function
$\tilde\eta:\cG^{-\infty}_{+,\sus}(M/B;E)\longrightarrow \bbC$ and the
discussion of multiplicativity shows that
\begin{equation}
\eta=\eta^+\circ q\circ\nu+\tilde\eta\circ\tq.
\label{perdet.177}\end{equation}
From the fundamental theorem of calculus,  
\begin{equation}
i\pi d\tilde\eta= \pi_{\infty}^*d\log\det
\label{perdet.179}\end{equation}
so we conclude from \eqref{perdet.177} that 
\begin{equation}
\tau=  e^{i\pi \eta} = e^{i\pi\eta^+\circ q\circ \pi}\;(\pi_{\infty}\tq)^*\det
\label{perdet.180}\end{equation}
defines the same cohomology class as the determinant on the classifying
group, \ie the first odd Chern class of the index bundle.
\end{proof}

\section{Adiabatic $\eta$}\label{Adeta}

We may further extend the discussion above by replacing the
once-product-suspended spaces by $(2n+1)$-times product-suspended spaces
using the isotropic quantization in $2n$ of the variables, as in
Theorem~\ref{fipomb2.16} applied to a decomposition
$\bbR^{2n+1}=\bbR\times\bbR^{2n}$ with the standard symplectic form used on
$\bbR^{2n}.$ Let $\cA[[\epsilon]]$ be the principal bundle of invertible
perturbations for the family $A$ with respect to the star product from
\eqref{iso.17}. 

\begin{proposition}\label{perdet.181} If
  $A\in\Psi^{m,m'}_{\psusn{2n+1}}(M/B;E,F)$ is a fully elliptic family and
  \eqref{perdet.50} is used, with the product interpreted as the
  parameter-dependent product of Theorem~\ref{fipomb2.16} for the
  symplectic form on $\bbR^{2n}$ then the resulting eta invariant on the
  bundle of smoothing pertubations has an asymptotic expansion as $\epsilon
  \downarrow 0$ which projects to
\begin{equation}
\eta_\epsilon :\cA[[\epsilon ]]\longrightarrow \epsilon ^{-n}\bbC[[\epsilon ]]
\label{perdet.182}\end{equation}
which has constant term the adibatic eta invariant 
\begin{equation}
\eta_{\adn{n}}:\cA[[\epsilon ]]\longrightarrow \bbC
\label{perdet.183}\end{equation}
which generates the first odd Chern class of the index bundle.
\end{proposition}

\begin{proof} This is essentially a notational extension of the results above.
\end{proof}

In particular \eqref{perdet.141} is a consequence of this result and Bott
periodicity. Namely, give an $2n$ product-suspended family we may always
choose a smoothing family, analogous to $Q_+$ in \eqref{perdet.171} which
is Schwartz in the second $2n-1$ variables and in the first is Schwartz at
$-\infty$ and of the form $Q_0+Q'$ with $Q'$ Schwartz at $+\infty$ and
$Q_0$ constant in the first variable (and Schwartz in the remainder). By
Bott periodicity, the even index of the family is the odd K-class on
$\bbR^{2n-1}\times B$ given by the product $(A(t)+Q_0)A(t)^{-1}$ for $t$
large. Then \eqref{perdet.141} follows by an elementary computation and the
proof of Lemma~\ref{detn.2} follows directly.

\appendix

\section{Symbols and products} \label{sp.0}

By choice of a quantization map, spaces of pseudodifferential operators on
a compact manifold can be identified, modulo smoothing operators, with the
appropriate spaces of symbols on the cotangent bundle as in
\eqref{fipomb2.7}. It is important to discuss, and carefully distinguish
between, several classes of such symbols and operators. To prepare for this
we describe here classes of product-type symbols for a pair of vector
spaces; subsequently this is extended to the case of vector bundles.

For a real vector space $V,$ the space of classical symbols of order $0$ on
$V$ is just $\CI(\com{V}),$ the space of smooth functions on the radial
compactification. In terms of any Euclidean metric on $V,$ $\rho
(v)=(1+|v|^2)^{-\frac12}$ is a defining function for the boundary of $\com{V}$
and the space of symbols of any complex order $z$ on $V$ is  
\begin{equation}
S^z(V)=\rho ^{-z}\CI(\com{V}).
\label{perdet.1}\end{equation}

If $W$ is a second real vector space then we may consider the radial
compactification $\com{V\times W}$ and corresponding symbol spaces
$S^z(V\times W).$ The natural projection $\pi_W:V\times W\longrightarrow W$
does not extend to a map from $\com{V\times W}$ to $\com{W}$ and
correspondingly classical symbols on $W$ do not generally lift to be
classical symbols on $V\times W.$ Rather $\com{V}\hookrightarrow\com{V\times
W}$ may be considered as an embedded submanifold, simply the closure (of
the preimage in $\com{V\times W})$ of $V\times\{0\}.$ On the other hand
there is certainly a smooth projection from $\com{V}\times\com{W}$ to
$\com{W};$ the smooth functions,
\begin{equation}
S^0(V;S^0(W))=S^0(W;S^0(V))=\CI(\com{V}\times\com{W}))
\label{perdet.2}\end{equation}
on this space are symbols on $V$ with values in the symbols on $W$ (or
vice-versa).

The main space we wish to consider here has some properties
between these two compactifications of $V\times W.$ Namely, in terms of
radial (real) blow-up, we set 
\begin{equation}
      \pcom{V}{V\times W}=[\com{V\times W};\pa(\com V\times\{0\})].
\label{perdet.3}\end{equation}
This manifold with corners has two boundary faces (unless one or both of the
factors is one-dimensional in which case either or both of the boundary
hypersurfaces may have two components). We use a subscript $V$ to refer to
the new boundary hypersurface produced by the blow-up in \eqref{perdet.3}.

\begin{lemma}\label{perdet.4} The projection $\pi_W:V\times W\longrightarrow W$
  extends to a smooth map
$$
\com{\pi}_W:\pcom{V}{V\times W}\longrightarrow \com{W}
$$
which is a fibration (with fibres which are manifolds with boundary) and in
terms of Euclidean metrics on $V$ and $W$ the functions
\begin{equation*}
\rho_V(v,w)=\left(\frac{1+|w|^2}{1+|v|^2+|w|^2}\right)^{\frac12}
\Mand\rho _r(v,w)=(1+|w|^2)^{-\frac12}
\label{perdet.5}\end{equation*}
extend from $V\times W$ to be smooth functions on $\pcom{V}{V\times W}$ and
are defining functions for the two boundary faces. 
\end{lemma}

\begin{proof} 
To check the first statement of the lemma, notice that the projection
$V\times W\to V$ has a smooth extension 
\begin{equation*}
p_{W}: \com{V\times W} \setminus \pa(\com{V}\times\{0\}) \to \com{W}
\end{equation*}
which is a fibration with typical fibre given by $V.$ Blowing up the
submanifold $\pa (\com{V}\times \{0\})$ in $\com{V\times W}$ exactly allows
us to extend $p_{W}$ to a fibration 
\begin{equation*}
\com{\pi}_W:\pcom{V}{V\times W}\longrightarrow \com{W}
  \end{equation*}
with typical fibre given by $\com{V}.$ Indeed, in $\com{V\times W}$ near
the submanifold $\pa(\com{V}\times\{0\})$, we can consider the generating
functions (\ie everywhere containing a coordinate system)
\begin{equation*}
\hat{v}= \frac{v}{|v|},\ \sigma_{V}= \frac{1}{(1+|v|^{2})^{\frac{1}{2}}},
\ \widetilde{w}=\frac{w}{(1+|v|^{2})^{\frac{1}{2}}}=\sigma_{V}w.
\end{equation*}
The blow up amounts to introducing polar coordinates
  \begin{equation*}
r_{V}= (\sigma_{V}^{2}+\widetilde{w}^{2})^{\frac{1}{2}},\ (\varphi,
\hat{\theta})=(\frac{\sigma_{V}}{r_{V}}, \frac{\widetilde{w}}{r_{V}})
\end{equation*}
so that the blow-down map is given locally by
  \begin{equation*}
\pcom{V}{V\times W}=[\com{V\times W};\pa(\com V\times\{0\})]
\ni (\hat{v}, r_{V}, \varphi, \hat{\theta})\longmapsto (\hat{v},
\sigma_{V}= r_{V}\varphi, \widetilde{w}= r_{V} \hat{\theta}).
\end{equation*}
In these polar coordinates, and for $r_{V}>0,$ the fibration $p_{W}$ is given by
\begin{equation}
p_{W}(\hat{v}, r_{V}, \varphi, \hat{\theta})=
\left(\frac{\varphi}{(\varphi^{2}+ |\hat{\theta}|^{2})^{\frac{1}{2}}},
\frac{\hat{\theta}}{(\varphi^{2}+ |\hat{\theta}|^{2})^{\frac{1}{2}}}\right)
\in \com{W}
 \label{ps.15}\end{equation}
where we have used the identification of $\com{W}$ with the upper
half-sphere which is the closure of the image
\begin{equation*}
W\ni w\longmapsto
(\frac{1}{(1+|w|^{2})^{\frac{1}{2}}},\frac{w}{(1+|w|^{2})^{\frac{1}{2}}})\in
\{(a,b)\in\bbR\times W; \quad a\ge 0, \;\; a^{2}+ |b|^{2}=1\}.
  \end{equation*} 
Thus, $p_{W}$ extends to $r_{V}=0$ to give the desired fibration.  
  
It follows from this that a defining function for the boundary of $\com{W}$
such as $(1+|w|^2)^{-\frac12}$ lifts from $\com{W}$ to be smooth and to
define the `old' boundary hypersurface, the one not produced by the blow
up.  Now $(1+|w|^2+|v|^2)^{\frac12}$ is a smooth boundary defining function
on $\com{V\times W}.$ It therefore lifts under the blow up in
\eqref{perdet.3} to be the product of defining functions for both boundary
hypersurfaces and so
\begin{equation*}
\rho_V(v,w)=\left(\frac{1+|w|^2}{1+|v|^2+|w|^2}\right)^{\frac12}
\end{equation*} 
is a boundary defining function for the new boundary produced by the blow-up.
\end{proof}

Now we define general spaces of `partial-product' symbols by 
\begin{equation}
S^{z,z'}(\pcom{V}{V\times W})=\rho_r^{z}\rho_V^{z'}\CI(\pcom{V}{V\times W}).
\label{perdet.6}\end{equation}
Directly from this definition,  
\begin{equation}
S^{z,z'}(\pcom{V}{V\times W})\cdot S^{\zeta,\zeta'}(\pcom{V}{V\times W})=
S^{z+\zeta,z'+\zeta'}(\pcom{V}{V\times W}).
\label{perdet.9}\end{equation}
Two of the `remainder' classes have simpler characterizations. Namely 
\begin{equation}
\begin{gathered}
S^{-\infty,z'}(\pcom{V}{V\times W})=\dCI(\com{W};S^{z'}(V))\\
S^{-\infty,-\infty}(\pcom{V}{V\times W})=\dCI(\com{V\times W})=\mathcal{S}(V\times W).
\end{gathered}
\label{perdet.10}\end{equation}

\section{Product Suspended operators}
\label{iso.0}

We can now introduce a generalization of the `suspended' algebra considered in
\cite{MR96h:58169} and in \cite{fipomb} (an algebra similar to the
suspended algebra was already introduced by Shubin in \cite{Shubin3}).  

The d-fold suspended pseudodifferential algebra on a compact manifold $X$
may be viewed as a space of smooth maps from $\bbR^d$ into $\Psi^k(X;E,F)$
in which the parameters (which we think of as the base variables for a
fibration) appear as `symbolic variables'. The inverse Fourier transform
identifies the suspended space
\begin{equation*}
\check{\Psi}_{\susn{d}}^k(X;E,F)\subset \Psi^k(\bbR^d\times X;E,F)
\label{fipomb2.18}\end{equation*}
directly, as is done in \cite{MR96h:58169}, with the
elements which are translation-invariant in $\bbR^d$ and have
convolution kernels vanishing rapidly at infinity, with all derivatives,
in these variables; this space may also be defined directly as in
\eqref{fipomb2.7}.

The subspace of smoothing operators is
\begin{equation*}
\Psi^{-\infty}_{\susn{d}}(X;E,F)=\cS(\bbR^d\times X^2;\Hom(E,F)\otimes\Omega _R)
\label{fipomb2.13}\end{equation*}
in terms of the Schwartz space. Then the finite-order operators may be
specified, up to smoothing terms, by Weyl quantization as
\begin{multline}
q_g:\rho ^{-k}\CI(\com{\bbR^d\times T^*X};\pi^*\hom(E,F))\ni a\longmapsto\\
(2\pi)^{-n}\int_{T^*X}\chi e^{iv(x,y)\cdot\xi}a(m(x,y),\zeta,\xi)d\xi dg
\in\Psi^k_{\susn{d}}(X;E,F)\label{fipomb2.14}\end{multline}
where the symbol space is compactified in the joint fibre 
$\bbR^d\times T^*_xX.$
The resulting full symbol sequence is as in \eqref{fipomb2.8} except that
the formal power series have coefficients on the sphere bundle of
$\bbR^p\times T^*X;$ the parameters do not affect the operators $B_j,$
acting on $T^*X,$ appearing in the product. 

If $A\in\Psi^{1}(X;E)$ 
is a first order pseudodifferential operator and $\tau$ is
the suspension variable for $\Psi^{1}_{\susn{1}}(X;E)$, then $A+i\tau$ is not in
general an element of $\Psi^{1}_{\susn{1}}(X;E).$ In fact, $A+i\tau\in
\Psi^{1}_{\susn{1}}(X;E)$ if and only if $A$ is a \emph{differential} operator.
Similarly, for $A\in\Psi^{1}(X;E,F)$, the operator 
\begin{equation}
\left(\begin{array}{cc}
   it+\tau & A^{*} \\
    A  & it-\tau
\end{array}\right)
\label{bp.1}\end{equation}
is in $\Psi^{1}_{\susn{2}}(X;E\oplus F)$ if and only if
$A$ is a differential operator.

This restriction to differential operators is unfortunate since the
operator $A+i\tau$ arises in the alternative definition of the eta
invariant as described in Appendix~\ref{eta.0}, while in
section~\ref{detn.0} the operator \eqref{bp.1} is used to implement Bott
periodicity for determinant line bundles.  For these reasons, and others,
we pass to the wider context of product-suspended operators.

We first need to enlarge the space of symbols as in section~\ref{sp.0}.  
Identifying $X$ with the zero section of $T^{*}X$, consider the 
blown-up space 
\begin{equation}
\pcom{X}{\bbR^{d}\times T^{*}X}=
[\com{\bbR^{d}\times T^{*}X};\pa\com{\bbR^{d}}\times X]
\label{sp.1}\end{equation}
where $\com{\bbR^{d}\times T^{*}X}$ is the radial compactification of 
$\bbR^{d}\times T^{*}X$ fibre by fibre and 
\begin{equation*}
\com{\bbR^{d}}\times X \subset \com{\bbR^{d}\times T^{*}X} 
\end{equation*}
is the closure of $\bbR^{d}\times X$ in $\com{\bbR^{d}\times T^{*}X}.$
In terms of a Riemannian metric $g$ and the Euclidean metric on $\bbR^{d}$.  
Lemma~\ref{perdet.4} generalizes directly to

\begin{lemma}
The projection $\bbR^{d}\times T^{*}X\to T^{*}X$ extends to a smooth map
\begin{equation*}
\com{\pi}_{T^*X}:\pcom{X}{\bbR^{d}\times T^{*}X}\longrightarrow\com{T^{*}X}
\end{equation*}
which is a fibration with typical fibre $\com{\bbR^{d}}$ and the smooth
functions 
\begin{equation*}
\bs(v,w)= \frac{(1+|w|^{2})^{\frac{1}{2}}}
{(1+|v|^{2}+|w|^{2})^{\frac{1}{2}}},\
\br(v,w)= (1+|w|^{2})^{-\frac{1}{2}},\ v\in \bbR^{d},\ w\in T^{*}X,
\end{equation*}
define the two boundary faces.  
\label{ps.2}\end{lemma}

\begin{proof}
This results from the invariance of the construction in \S\ref{sp.0} under those
linear transformations of $V\times W$ which leave $V$ invariant, so
Lemma~\ref{perdet.4} extends to the case of a vector bundle.
\end{proof}

For $z,z'\in \bbC,$ the space of (partially) product-type symbols
with values in a vector bundle over $X$ is then
\begin{equation}
\cS^{z,z'}(\pcom{X}{\bbR^{d}\times T^{*}X};U)=
\br^{-z}\bs^{-z'}\CI(\pcom{X}{\bbR^{d}\times T^{*}X};U).
\label{ps.3}\end{equation}
On $\com{\bbR^{d}}\times X\times X$, consider the boundary defining
function $\bt(\tau)=(1+|\tau|^{2})^{-\frac{1}{2}}.$ Let $E$ and $F$ be
smooth complex vector bundles on $X.$ For $z'\in\bbC$ set
\begin{equation}
\Psi^{-\infty,z'}_{\psusn{d}}(X;E,F)=
\bt^{-z'}\CI(\com{\bbR^{d}}\times X\times X; \Hom(E,F)\otimes \Omega_{R}X),
\label{ps.4}\end{equation}
where $\Omega_{R}X= \pi_{3}^{*}\Omega X,$ $\pi_{3}$ being the projection
on the third factor, and $\Omega X$ being the bundle of densities on $X.$
This is the space of smoothing operators (defined as usual through their
kernels) on $X$ depending symbolically on $d$ parameters; which we identify
as the product-suspended operators of order $-\infty$ on $X.$

\begin{definition} The general spaces of product $d$-suspended
pseudodifferential operators of order $k,k'\in\bbZ$ acting from 
$\cS(\bbR^{d}\times X;E)$ to $\cS(\bbR^{d}\times X;F)$ is
\begin{equation*}
  \Psi_{\psusn{d}}^{k,k'}(X;E,F)=
   q_{g}(\cS^{k,k'}(\pcom{X}{\bbR^{d}\times T^{*}X};\hom(E,F))) +
  \Psi_{\psusn{d}}^{-\infty,k'}(X;E,F)
\end{equation*} 
where $q_{g}$ is the Weyl quantization \eqref{fipomb2.14} applied to these more
general symbol spaces.
\label{ps.6}\end{definition}
\noindent We limit attention to integral orders here only because it is all
that is needed.

Pseudodifferential operators are included in the product-suspended operators
\begin{equation*}
\Psi^{k}(X;E,F) \subset \Psi^{k,0}_{\psusn{d}}(X;E,F),
\end{equation*}
being independent of the parameters. For integers $l\le k,$ $l'\le k'$,
there are inclusions
\begin{equation*}
     \Psi_{\psusn{d}}^{l,l'}(X;E,F)\subset \Psi_{\psusn{d}}^{k,k'}(X;E,F).    
\end{equation*}
Furthermore, as we will see below in Theorem~\ref{fipomb2.16},
product $d$-suspended operators compose in the expected way
\begin{equation*}
       \Psi_{\psusn{d}}^{k,k'}(X;E,F)\circ \Psi_{\psusn{d}}^{l,l'}(X;G,E)
\subset \Psi_{\psusn{d}}^{k+l,k'+l'}(X;G,F).
\end{equation*}
  
Suspended operators are particular instances of  product suspended operators,
\begin{equation*}
\Psi_{\susn{d}}^{k}(X;E,F)\subset  \Psi_{\psusn{d}}^{k,k}(X;E,F),\ k\in\bbZ, 
\end{equation*}
and 
\begin{equation*}
\Psi_{\susn{d}}^{-\infty}(X;E,F)=\Psi_{\psusn{d}}^{-\infty,-\infty}(X;E,F).
\end{equation*}
\label{prsu.4}

Product $d$-suspended pseudodifferential operators are intimately
related with the algebra of product-type operators introduced in
\cite{faficu}.  More precisely, consider the projection
\begin{equation}
\phi: \bbR^{d}\times X \to \bbR^{d} 
\label{pspp.6}\end{equation} 
as a fibration. If $E$ and $F$ are smooth complex vector bundles on $X,$
then as discussed in \cite{faficu}, to such a fibration one can associate
the space of product-type pseudodifferential operators of order $(k,k')$
\begin{equation*}
      \Psi_{\pt}^{k,k'}(\bbR^{d}\times X;E,F)
\end{equation*}
acting from $\CI_{c}(\bbR^{d}\times X;E)$ to $\CI(\bbR^{d}\times X;F)$.
Given $\tau\in \bbR^{d}$, let 
\begin{equation*}
T_{\tau}:\bbR^{d}\times X \to\bbR^{d}\times X
\end{equation*}
denote the translation in the first factor $T_{\tau}(t,x)= (t-\tau, x)$. 
We can consider the product-type psuedodifferential operators which are
translation-invariant in the Euclidean variable, that is, satisfying
\begin{equation}
T_{\tau}^{*}(Af)= AT^{*}_{\tau}f, \forall\; \tau \in \bbR^{d}, \,
   f\in \CI_{c}(\bbR^{d}\times X;E).
\label{ps.7}\end{equation}
In terms of the Schwartz kernel $K_{A}$ of $A$, this means that $K_{A}$ acts
by convolution in the first factor
\begin{equation*}
   Af(x,t)= \int_{\bbR^{d}}\int_{X} K_{A}(t-s, x,x') f(x',s) ds
\end{equation*} 
where $K_{A}$ is a density in the $x'$ variable. Now one can ask in
addition that this convolution kernel decay to all orders at infinity
\begin{equation}
  K_{A}\in \mathcal{C}^{-\infty}_{c}(\bbR^{d}\times X^{2};
   \Hom(E,F)\otimes \Omega_{R}X) + \cS(\bbR^{d}\times X^{2};
\Hom(E,F)\otimes \Omega_{R}X).
\label{ps.8}\end{equation}
This leads to following characterization of product $d$-suspended operators.

\begin{lemma}\label{ps.9} Fourier transformation in the suspension variables
\begin{equation*}
         (\hat{A}(\tau)f)(x)= \int_{X}\int_{\bbR^{d}} e^{-it\tau} 
                K_{A}(t,x,x')f(x')dt, \quad \tau\in \bbR^{d}
\end{equation*}
is an isomorphism of the space of translation-invariant product-type
pseudodifferential operators satisfying \eqref{ps.8} onto the $d$-parameter
product-suspended pseudodifferential operators; it preserves products.
\label{ps.10}\end{lemma}

\begin{proof}
Modulo small changes of notation, this is the same as for suspended operators.
\end{proof}

One advantage of the alternative definition through Lemma~\ref{ps.9} is
that the Fredholm theory for product $d$-suspended operators follows almost
immediatly from the corresponding Fredholm theory for product-type
operators.  Indeed, the principal symbol map and the base family map for
product-type operators gives via the inclusion (using the inverse Fourier
transform) $\check{\Psi}_{\psusn{d}}(X;E,F)\subset
\Psi_{\pt}^{k,k'}(\bbR^{d}\times X;E,F)$ a corresponding symbol map and
base family map for product $d$-suspended operators.  For the convenience
of the reader, we will define these directly without refering to
product-type operators.

Of the two boundary faces of $\pcom{X}{\bbR^p\times T^*X},$ the `old'
boundary, or really its blow-up,
\begin{equation*}
       B_{\sigma}=[S(\bbR^{d}\times T^{*}X); S(\bbR^{d})\times X ]
\end{equation*}
with $X$ being the zero section of $T^*X,$ carries the replacement for the
usual principal symbol. In terms of a quantization map as above, this is given
by the restriction of the full symbol $a\in \cS^{m,m'}(\bbR^{d}\ltimes
T^{*}X; E,F)$ of an operator $A= q_{g}(a)$ to this boundary face,
\begin{equation}
\sigma_{m,m'}:\Psi^{m,m'}_{\psusn{d}}(X;E,F)\longrightarrow
\cS^{m,m'}_{\psusn{d}}(X;E,F)
\label{ps.12}\end{equation}
with 
\begin{equation*}
  \cS^{m,m'}_{\psusn{d}}(X;E,F)= \CI(B_{\sigma};
\hom(E,F)\otimes N^{-m}\otimes N^{-m'}_{\ff}) 
\end{equation*}
where $N$ is the normal bundle to $B_{\sigma}$ and
and $N_{\ff}$ is the normal bundle of the `new' boundary, which is
canonically identified with the normal bundle to the boundary of
$B_{\sigma}.$ Both are trivial bundles. This corresponds to the
multiplicative short exact sequence
\begin{equation}
0\longrightarrow 
\Psi_{\psusn{d}}^{m-1,m'}(X;E,F)\longrightarrow 
\Psi_{\psusn{d}}^{m,m'}(X;E,F)
\overset{\sigma_{m,m'}}{\longrightarrow} \cS^{m,m'}_{\psusn{d}}(X;E,F)\to 0. 
\label{ps.13}\end{equation}

A product $d$-suspended operator $A\in \Psi_{\psusn{d}}^{k,k'}(X;E,F)$ is
 \emph{elliptic} if its principal symbol $\sigma_{m,m'}(A)$ is invertible.

Ellipticity alone does imply that the family is Fredholm for each value of
the parameter but, as for product-type operators, is does not suffice to
allow the construction of a parametrix modulo Schwartz-smoothing
errors. There is a second symbol map which takes into account the behavior of
the operator for large values of the suspension parameters.

Let $B_{\sus}\subset\pcom{X}{\bbR^{d}\times T^{*}X}$ denote the `new' boundary,
which is the `front face' produced by the blow up. The
fibration of Lemma~\ref{ps.2} gives a canonical identification of $B_{\sus}$ with
$S(\bbR^{d})\times\com{T^{*}X}.$ Thus, the restriction map (using a
boundary defining function $\rho _{\sus}$ for $B_{\sus})$ becomes
\begin{multline}
R:\cS^{k,k'}(\bbR^{d}\ltimes X; \hom(E,F)\ni a \longmapsto\\
\bs^{m'}a\big|_{B_{\sus}}\in
\CI(S(\bbR^{d});\cS^{k,k'}(\com{T^{*}X}; \hom(E,F)))
\label{ps.17}\end{multline}

Given any element $A=q_{g}(a_{1})+A_{2}\in
\Psi_{\psusn{d}}^{k,k'}(X;E,F)$ with $a_{1}\in \cS^{k,k'}(\bbR^{d}\ltimes
X; \hom(E,F)$ and  $a_{2}\in\Psi_{\psusn{d}}^{-\infty,k'}(X;E,F),$ the 
\emph{base family} is defined by
\begin{equation}
L(A)= q_{g}(R(a_{1}))+\rho_{\sus}^{k'}A_{2})\big|_{B_{\sus}}
\in\CI(S(\bbR^{d});\Psi^{m}(X;E,F))
\label{ps.18}
\end{equation}

\begin{proposition} The base family \eqref{ps.18} is independent of
  choices and corresponds to the multiplicative short exact sequence
\begin{equation}
0\longrightarrow
\Psi^{k,k'-1}_{\psusn{d}}(X;E,F)\longrightarrow
\Psi^{k,k'}_{\psusn{d}}(X;E,F)
\overset{L}{\longrightarrow}
\CI(S(\bbR^{d});\Psi^{k}(X;E,F))\longrightarrow
0,
\label{ps.19}\end{equation}
so
\begin{equation*}
L(A\circ B)= L(A)\circ L(B),\ A\in \Psi_{\psusn{d}}^{m,m'}(X;E,F),\
B\in \Psi_{\psusn{d}}^{k,k'}(X;G,E).
\end{equation*}
\end{proposition}\label{ps.23.03.2006}

\begin{proof}
The fact that there is a short exact sequence is essentially by definition of 
$L.$ The fact that $L$ is a homomorphism follows by very simple `oscillatory
testing'. Namely, if $u\in\CI(X;E)$ and
$A\in\Psi^{k,k'}_{\psusn{p}}(X;E,F)$ then 
\begin{equation}
Au\in\rho^{-k'}_{\tau}\CI(\com{\bbR^p}\times X;F)\Mand L(A)u=\rho
_{\tau}^{k'}Au\big|_{\pa{\com{\bbR^p}}}\in\CI(\bbS^{p-1}\times X;F).
\label{perdet.162}\end{equation}
\end{proof}

\begin{definition}
The \emph{joint symbol} $J(A)$ of an operator $A\in 
\Psi^{k,k'}_{\psusn{d}}(X;E,F)$ is the combination of its principal symbol and
its base family
\begin{equation*}
J(A)= (\sigma(A), L(A))\Mwhere\sigma(L(A))= \sigma(A)\big|_{B_{\sigma }}.
\end{equation*}
An operator $A$ is said to be \emph{fully elliptic} if its joint symbol is
invertible.
\label{js.1}\end{definition}

The important feature that motivates the introduction of product-suspended 
operators (as opposed to suspended operators) is the following lemma.

\begin{lemma} If $A\in \Psi^{1}(X;E)$ then the one-parameter family
$\tau\longmapsto A+i\tau\in\Psi^{1,1}_{\psus(1)}(X;E)$ and if $B\in
\Psi^{1}(X;E,F)$, then the two-parameter family
\begin{equation*}
(t,\tau)\longmapsto\hat{B}(t,\tau)= \begin{pmatrix}
it+\tau & B^{*} \\
B & it-\tau 
\end{pmatrix}\in\Psi^{1,1}_{\psus(2)}(X;E\oplus F).
\end{equation*}
Moreover if $A$ is self-adjoint and elliptic (respectively 
$B$ is elliptic) then $A+i\tau$ (respectively $\hat{B}$) is fully elliptic. 
\label{prsu.5}\end{lemma}
\noindent In fact, it suffices that all the eigenvalues of the symbol of
$A$ have a nonvanishing real part for $A+i\tau$ to be fully elliptic.
 
\begin{proof}
Fix a quantization $q_{g}.$ In the first case $a\in
\rho^{-1}\CI(\com{T^{*}X}; \pi^{*}\hom(E))$ exists such that
$(A-q_{g}(a))\in\Psi^{-\infty}(X;E).$ Then 
\begin{equation*}
a+i\tau \in \cS^{1,1}(\bbR\ltimes T^{*}X;E)\Mand
   A+i\tau -q_{g}(a+i\tau)\in \Psi_{\psusn{1}}^{-\infty,1}(X;E),
\end{equation*}
which shows that $A+i\tau\in \Psi^{1,1}_{\psusn{1}}(X;E).$ The symbol of
$A+i\tau$ is invertible if $\sigma (A)$ has no eigenvalues in $i\bbR$ and
its base family is $\pm i\Id$ at the two components of
$\pa\com{\bbR_{\tau}}.$ Thus $A+i\tau$ is fully elliptic.

In the second case, choose $b\in \rho^{-1}\CI(\com{T^{*}X};
\pi^{*}\hom(E,F))$ such that
\begin{equation*}
            B-q_{g}(b)\in \Psi^{-\infty}(X;E).
\end{equation*}
Then 
\begin{equation*}
   \hat{b}= \begin{pmatrix}it+\tau & b^{*} \\
b & it-\tau 
\end{pmatrix} \in \cS^{1,1}(\bbR^{2}\ltimes T^{*}X;E,F)
\end{equation*}
and $\hat{B}- q_{g}(\hat{b})\in \Psi^{-\infty,1}_{\psusn{2}}(X;E,F),$
which shows that $\hat{B}\in\Psi^{1,1}_{\psusn{1}}(X;E,F).$ To see that
$\hat{B}$ is fully elliptic when $B$ is elliptic, 
consider the invertible operator
\begin{equation*}
 Q= \hat{B}^{*}\hat{B}+1=\begin{pmatrix}
B^{*}B+t^{2}+\tau^{2}+1 & 0 \\
0 &  BB^{*}+t^{2}+\tau^{2}+1
\end{pmatrix} \in \Psi^{2,2}_{\psusn{2}}(X;E\oplus F).
\end{equation*}
Then
\begin{equation*}
     (Q^{-1}\hat{B}^{*})\hat{B} -\Id_{E\oplus F}=
          -Q^{-1}\in \Psi^{-2,-2}_{\psusn{2}}(X;E\oplus F),
\end{equation*}
so that $J(\hat{B})^{-1}= J(Q^{-1}\hat{B}^{*})$ exists, which shows that
$\hat{B}$ is fully elliptic.
\end{proof}

\section{Mixed isotropic operators}\label{wq.0}

Next we proceed to the `parameter quantization' of these spaces of product
suspended operators. That is we introduce a new product depending on the
choice of an antisymmetric form on $\bbR^{p}.$ These products are used
above in the identification of the determinant bundle, as constructed in
the product $2n$-suspended case, with the determinant bundle as introduced
by Quillen. To do so we use an adiabatic limit, with a parameter which
passes from the quantized to the unquantized case discussed above; for the
isotropic algebra itself such degenerations are treated in \cite{HHH2} and
as shown there implements Bott periodicity. So, to introduce these spaces
we simply combine \eqref{fipomb2.7} and its Euclidean analogue
\eqref{fipomb2.11}.  Note that the quantization map will be global in the
Euclidean variables but can only be local near the diagonal in the
manifold. In defining these spaces we use the formula for the action of an
operator by Weyl quantization in \eqref{perdet.161}.

\begin{proposition}\label{iso.1} Let $X$ be a compact manifold $E$ and $F$
complex bundles over $X$ then for any $p\in\bbN$ combining
\eqref{perdet.161} with the operator product gives a smooth family of associative
products 
\begin{equation}
\Lambda^2(\bbR^n)\times\Psi_{\psus(2n)}^{m_1,m'_1}(X;F,G)\times
\Psi_{\psus(2n)}^{m_2,m'_2}(X;E,F)\longrightarrow
\Psi_{\psus(2n)}^{m_1+m_2,m'_1+m'_2}(X;E,G).
\label{fipomb2.19}\end{equation}
\end{proposition}

This follows by combining essentially standard treatments of
the composition of pseudodifferential operators with those of the
`isotropic' operators on $\bbR^n.$

We are especially interested in the `adiabatic limit' where the general
$\omega$ is replaced by $\epsilon \omega$ for a fixed antisymmetric
form. The cases which occur above are where $p$ is even and $\omega$ is
non-degenerate, or where $p$ is odd and $\omega$ has maximal rank. In this
case we state the corresponding corollary of the result above (see also
\cite{Hormander5} and \cite{HHH2}).
   
\begin{theorem}\label{fipomb2.16} For any fixed antisymmetric form on
  $\bbR^p,$ the composition \eqref{fipomb2.19} induces a smooth 1-parameter
  family of quantized products
 \begin{equation}
[0,1]_{\epsilon}\times 
\Psi^{k,k'}_{\psus(p)}(X;F,G)\times\Psi^{l,l'}_{\psus(p)}(X;E,F)
\longrightarrow \Psi^{k+l,k'+l'}_{\psus(p)}(X;E,G)
\label{iso.8}\end{equation}
and as $\epsilon \downarrow0$ there is a Taylor series expansion 
\begin{equation}
A \circ_{\epsilon} B(u)\sim
\sum_{k=0}^{\infty}\frac{(-i\epsilon)^{k}}{2^{k}k!}
\omega(D_{v},D_{w})^{k}
A(v)B(w) \big|_{v=w=u}'
\label{iso.17}\end{equation}
in particular, when $\epsilon=0$ the product reduces to the usual
parameterized product of suspended operators.
\end{theorem}

\providecommand{\bysame}{\leavevmode\hbox to3em{\hrulefill}\thinspace}
\providecommand{\MR}{\relax\ifhmode\unskip\space\fi MR }
\providecommand{\MRhref}[2]{%
  \href{http://www.ams.org/mathscinet-getitem?mr=#1}{#2}
}
\providecommand{\href}[2]{#2}

\end{document}